\documentclass[a4paper,12pt,reqno]{amsart}

\usepackage{amssymb, amsmath, amsthm, mathrsfs, braket, xspace}
\usepackage[margin=1.1in]{geometry} 

\numberwithin{equation}{section}
\usepackage[colorlinks=true]{hyperref}

\usepackage{mathtools}
\mathtoolsset{showonlyrefs=true} 
\usepackage{enumitem}
\setlist[enumerate]{label=$(\mathrm{\arabic*})$,leftmargin=*}

\usepackage[all,2cell]{xy}
\objectmargin+{1mm}
\labelmargin+{0.8mm}
\SelectTips{cm}{12}
\usepackage{tikz}
\usepackage{float}
\usepackage{aliascnt}
\newtheorem{thm}{Theorem}[section]

\newaliascnt{cor}{thm}
\newtheorem{cor}[cor]{Corollary}
\aliascntresetthe{cor}

\newaliascnt{lem}{thm}
\newtheorem{lem}[lem]{Lemma}
\aliascntresetthe{lem}

\newaliascnt{prop}{thm}
\newtheorem{prop}[prop]{Proposition}
\aliascntresetthe{prop}

\newaliascnt{conj}{thm}

\aliascntresetthe{conj}

\theoremstyle{definition}
\newaliascnt{dfn}{thm}
\newtheorem{dfn}[dfn]{Definition}
\aliascntresetthe{dfn}

\newaliascnt{rem}{thm}

\aliascntresetthe{rem}

\newaliascnt{prob}{thm}

\aliascntresetthe{prob}

\newaliascnt{ex}{thm}

\aliascntresetthe{ex}

\newcommand{\Z}{\mathbb{Z}}

\newcommand{\bA}{\mathbb{A}}

\newcommand{\bP}{\mathbb{P}}

\newcommand{\cF}{\mathcal{F}}

\DeclareMathOperator{\Coker}{Coker}

\DeclareMathOperator{\dlog}{dlog}

\DeclareMathOperator{\Hom}{Hom}
\DeclareMathOperator{\img}{Im}

\DeclareMathOperator{\Ker}{Ker}
\DeclareMathOperator{\res}{res}
\DeclareMathOperator{\Spec}{Spec}
\DeclareMathOperator{\Tr}{Tr}

\newcommand{\Gm}{\mathbb{G}_{m}}
\newcommand{\Ga}{\mathbb{G}_{a}}
\newcommand{\Id}{\mathrm{id}}

\newcommand{\bcube}{{\overline{\square}}}

\renewcommand{\d}{\mathrm{d}}

\newcommand{\black}{\color{black}}

\newcommand{\RSC}{{\operatorname{\mathbf{RSC}}}}
\newcommand{\CI}{\operatorname{\mathbf{CI}}}
\newcommand{\MPST}{\operatorname{\mathbf{MPST}}}
\newcommand{\ulMPST}{\operatorname{\mathbf{{\underline{M}}PST}}}
\newcommand{\PST}{\operatorname{\mathbf{PST}}}
\newcommand{\HI}{\operatorname{\mathbf{HI}}}
\newcommand{\Sch}{\operatorname{\mathbf{Sch}}}
\newcommand{\Sm}{\operatorname{\mathbf{Sm}}}
\newcommand{\Cork}{\mathbf{Cor}}
\newcommand{\MCor}{\operatorname{\mathbf{MCor}}}
\newcommand{\ulMCor}{\operatorname{\mathbf{{\underline{M}}Cor}}}
\newcommand{\Acal}{\mathcal{A}}
\newcommand{\Dcal}{\mathcal{D}}
\newcommand{\Pcal}{\mathcal{P}}
\newcommand{\Ztr}{\Z_{\mathrm{tr}}}

\newcommand{\GaMod}{\mathbb{G}_{a}^{M}}
\newcommand{\GmMod}{\mathbb G_m^M}
\newcommand{\GaSharp}{\mathbb{G}_{a}^{\#}}
\newcommand{\GmSharp}{\mathbb G_m^{\#}}
\newcommand{\Nis}{\mathrm{Nis}}

\newcommand{\Gapair}{\mathcal{G}_a}
\newcommand{\Gmpair}{\mathcal{G}_m}

\newcommand{\va}{\mathbf{a}}
\newcommand{\vb}{\mathbf{b}}
\newcommand{\vc}{\mathbf{c}}
\newcommand{\vd}{\mathbf{d}}

\title{Tensor products of additive groups as reciprocity sheaves}

\author[T. Hiranouchi]{Toshiro Hiranouchi}\address[T. Hiranouchi]{
Department of Basic Sciences, Graduate School of Engineering,
Kyushu Institute of Technology,
1-1 Sensui-cho, Tobata-ku, Kitakyushu-shi,
Fukuoka, 804-8550 JAPAN}
\email{hira@mns.kyutech.ac.jp}

\author[R. Sugiyama]{Rin Sugiyama} \address[R. Sugiyama]{
Department of Mathematics, Physics and Computer Science,
Japan Women's University, 2-8-1 Mejirodai, Bunkyo-ku, Tokyo, 112-8681 JAPAN
}\email{sugiyamar@fc.jwu.ac.jp}


\begin{document}

\begin{abstract}
We compute the Nisnevich reciprocity sheaves associated with
tensor products of the cube-invariant objects for the additive
and multiplicative groups, with at least one additive factor.
Over a perfect field of characteristic different from two, these sheaves 
decompose into direct sums of absolute K\"ahler differential sheaves with binomial
multiplicities.  
As an application, we compute the associated reciprocity sheaves
of the Hochschild homology of the additive object.
\end{abstract}

\maketitle

\section{Introduction}\label{sec:intro}
In the category $\HI$ of homotopy invariant sheaves with transfers, by the work of many authors, 
the tensor product of the multiplicative group $\Gm$ 
is closely related to motivic cohomology. 
More precisely,
for $q\geq1$, there are isomorphisms
\begin{align}\label{eq:Gmtensor}
    \Gm^{\overset{\HI}{\otimes}q}\simeq \mathcal{H}^q(\Z(q))\simeq \mathcal{K}_q^M
\end{align}
where $\Z(q)$ is Voevodsky's motivic complex and $\mathcal{K}^M_q$ is the improved Milnor $K$-sheaf (\cite{Ker10}).
The tensor product in $\HI$ 
has a useful
extension to reciprocity sheaves through modulus presheaves with transfers.
In \cite[Sections~4--5]{RSY22}, R\"ulling, Sugiyama and Yamazaki developed this construction and computed a number of basic examples. 
In particular,
\cite[Proposition~5.6]{RSY22} recovers \eqref{eq:Gmtensor} using a cube-invariant object $\GmSharp$ whose associated reciprocity sheaf is $\Gm$;
\[
h_{0,\Nis}\left((\GmSharp)^{\otimes q}
 \right)
 \simeq
 \mathcal{K}^M_q.
\]
Furthermore, \cite[Theorem~5.14 and Theorem~5.15]{RSY22} computed the tensor product of two cube-invariant objects 
$\GaSharp$ and 
related 
the additive group $\Ga$ and the sheaf $\Omega^1$ of absolute K\"ahler differentials of degree one: 
\[
 h_{0,\Nis}(\GaSharp\otimes\GaSharp)
 \simeq \Ga\oplus\Omega^1.
\]
The purpose of this paper is to extend this computation to arbitrary
tensor powers of $\GaSharp$ and their mixed tensor products with $\GmSharp$. 
We show that the associated Nisnevich
reciprocity sheaves decompose into direct sums of
absolute K\"ahler differential sheaves.

To state our results more precisely, 
let $k$ be a perfect field of characteristic $\neq 2$. 
Consider the modulus pairs  
$\Gapair=(\bP^1,2\infty)$, and $\Gmpair=(\bP^1,(0)+(\infty))$ defined over $k$.
Following \cite[Section~5.1]{RSY22}, let
$\GaSharp$ and $\GmSharp$ be the reduced cube-invariant objects associated with
$\Gapair$ and $\Gmpair$.  We write $\GaMod$ and $\GmMod$ for the pointed
modulus presheaves obtained from these pairs by the zero section and
the unit section, respectively (for the precise definitions, see \autoref{sec:prelim}). 
For a modulus presheaf or a cube-invariant object
$H$, $h_{0,\Nis}(H)$ denotes its associated Nisnevich reciprocity sheaf.  
The  associated reciprocity sheaves of the cube-invariant objects $\GaSharp$ and $\GmSharp$ are 
the additive group $\Ga$ and the multiplicative group $\Gm$, respectively; 
\[
h_{0,\Nis}(\GaSharp) \simeq \Ga,\quad h_{0,\Nis}(\GmSharp) \simeq \Gm
\]
(\cite[Section~5.1]{RSY22}). 
Throughout this paper,
$\otimes$ denotes the tensor product induced by the symmetric monoidal
structure (see \autoref{sec:prelim} for the precise definitions of these objects).  
Our main result is the following.

\begin{thm}[{\autoref{thm:decomp}}]\label{thm:main}
Let $k$ be a perfect field of characteristic $\neq 2$, and let 
$n,m\geq 0$.
There are canonical isomorphisms of Nisnevich reciprocity sheaves
\begin{equation}\label{eq:maini}
 \begin{split}
 h_{0,\Nis}\left(
 (\GaMod)^{\otimes (n+1)}\otimes(\GmMod)^{\otimes m}
 \right)
 &\simeq
 h_{0,\Nis}\left(
 (\GaSharp)^{\otimes (n+1)}\otimes(\GmSharp)^{\otimes m}
 \right)\\
 &\simeq
 \bigoplus_{q=0}^{n}
 (\Omega^{m+q})^{\oplus\binom{n}{q}}.
 \end{split}
\end{equation}

\end{thm}

The following low-degree cases illustrate \autoref{thm:main}:
\begin{align*}
 h_{0,\Nis}((\GaSharp)^{\otimes 2})&\simeq \Ga\oplus\Omega^1,\\
 h_{0,\Nis}((\GaSharp)^{\otimes 3})&\simeq
 \Ga\oplus(\Omega^1)^{\oplus 2}\oplus\Omega^2.
\end{align*}
The first nontrivial case is exactly the computation of $h_{0,\Nis}((\GaSharp)^{\otimes 2})$ in
\cite[Theorem~5.14]{RSY22}.

As an application, we study in \autoref{sec:HH} the Hochschild homology of
$\GaSharp$.  Multiplication of affine coordinates makes $\GaSharp$ an
associative algebra object in the category $\CI$ of cube-invariant objects. 
If $\mathrm{HH}_n^{(m)}(\GaSharp)$ denotes
the $n$-th homology object of the Hochschild complex
\[
C_n^{(m)}(\GaSharp) :=
 (\GaSharp)^{\otimes (n+1)}\otimes(\GmSharp)^{\otimes m}, 
\]
we obtain the following result.

\begin{thm}[{\autoref{thm:HH}}]\label{thm:HHin}
Let $k$ be a perfect field of characteristic $p\neq 2$. 
For every $n,m\geq0$, there is a canonical isomorphism of Nisnevich
reciprocity sheaves
\begin{equation}\label{eq:HH-intro}
 h_{0,\Nis}\bigl(\mathrm{HH}_n^{(m)}(\GaSharp)\bigr)\simeq\Omega^{n+m}.
\end{equation}
\end{thm}

Thus, after passage to Nisnevich reciprocity sheaves, the Hochschild
homology of $\GaSharp$ recovers the sheaves of absolute K\"ahler
differentials.

\autoref{app:diff} derives from \cite[Lemma~4.2]{BK86} an
$F$-linear presentation of absolute K\"ahler differentials over a field
$F$ and compares it with \cite[Lemma~4.1]{BE03}.  This $F$-linear form is
used in the construction of the inverse maps below.

\subsection*{Notation}
Throughout the paper, $k$ is a perfect field of characteristic $p \neq 2$.  For a
field $F/k$, we write $\Omega_F^q=\Omega^q_{F/\Z}$ for the absolute K\"ahler
differentials, and $\Omega^q$ for the corresponding Nisnevich sheaf.  For a
finite extension $E/F$, the usual trace and the differential trace are both
denoted by $\Tr_{E/F}$; the source makes the meaning clear.

\subsection*{Acknowledgements}
The first author was supported by JSPS KAKENHI Grant Number 24K06672.
The second author was supported by JSPS KAKENHI Grant Number 21K03188.

\section{Preliminaries}\label{sec:prelim}
The theory of modulus presheaves with transfers and reciprocity sheaves has been developed by a series of works of Kahn-Miyazaki-Saito-Yamazaki \cite{KSY16, KSY22, KMSY1, KMSY2}.
We recall only the part of the theory of modulus presheaves that is needed
below.  A convenient short introduction with the same conventions is given
in \cite[Section~5, especially Definitions~5.1--5.7]{HS25}; for the full
tensor construction we use \cite[Sections~2--4]{RSY22}.

\subsection{Modulus presheaves and cube-invariance}
Let $\Sch$ be the category of separated schemes of finite type over $k$
and let $\Sm\subset\Sch$ be its full subcategory of smooth schemes.  For
$X\in\Sm$ and $Y\in\Sch$, a \emph{prime correspondence} from $X$ to $Y$
is an integral closed subscheme $\Gamma\subset X\times Y$ which is finite
and surjective over a connected component of $X$.  Following
\cite[Section~2.1]{RSY22}, we use the term prime correspondence; the same
objects are called elementary correspondences in
\cite[Definition~1.1]{MVW06}.  For $X,Y\in\Sm$, the group
$\Cork(X,Y)$ is the free abelian group on the prime correspondences from
$X$ to $Y$.  We write $\PST$ for the category of presheaves with transfers
on $\Cork$.  For $X\in\Sm$, put
\begin{equation}\label{eq:ZtrSm}
 \Ztr(X):=\Cork(-,X)\in\PST.
\end{equation}
See \cite[Sections~1.2 and~2.1]{RSY22}.

A \emph{modulus pair} is a pair $M=(\overline M,M^\infty)$, where
$\overline M\in\Sch$, $M^\infty$ is an effective Cartier divisor on
$\overline M$, and
$M^\circ:=\overline M\smallsetminus |M^\infty|$ is smooth.  It is called
proper if $\overline M$ is proper over $k$.  Let
$M=(\overline M,M^\infty)$ and $N=(\overline N,N^\infty)$ be modulus
pairs and let $\Gamma\in\Cork(M^\circ,N^\circ)$ be a prime
correspondence.  Write $\overline\Gamma\subset\overline M\times\overline N$
for its closure and $\overline\Gamma^N\to\overline\Gamma$ for the
normalization.  The correspondence $\Gamma$ is \emph{admissible} if
\begin{equation}\label{eq:admissdef}
 M^\infty|_{\overline\Gamma^N}
 \geq
 N^\infty|_{\overline\Gamma^N},
\end{equation}
where both divisors are pulled back from $\overline M$ and $\overline N$,
respectively.  
It is \emph{left proper} if $\overline\Gamma$ is proper
over $\overline M$.  The group $\ulMCor(M,N)$ is the subgroup of
$\Cork(M^\circ,N^\circ)$ generated by the admissible left proper prime
correspondences.  The category $\ulMCor$ has all modulus pairs
as objects and these groups as morphisms; see \cite[Section~2.2]{RSY22}.
We write
\[
 \MCor\subset \ulMCor
\]
for the full subcategory of proper modulus pairs, as in
\cite[Section~2.11]{RSY22}.  Accordingly,
$\ulMPST$ denotes the category of additive presheaves on
$\ulMCor$, whereas $\MPST$ denotes the category of additive
presheaves on $\MCor$, called modulus presheaves with transfers.

If $\overline M$ is normal and a morphism
$g^\circ:M^\circ\to N^\circ$ extends to a morphism
$g:\overline M\to\overline N$, then the admissibility of its graph is
simply the divisor inequality
\begin{equation}\label{eq:admissmap}
 M^\infty\geq g^*N^\infty.
\end{equation}
This is the form used below.  For example, for a regular proper curve $C/F$ over a field $F/k$
and $a\in F(C)$, the morphism
$(C,D)\to(\bP^1,2\infty)$ induced by $a$ is admissible whenever 
$D\geq 2\mathrm{div}_\infty(a)$ where $\mathrm{div}_\infty(a)$ is the pole divisor of $a$.

For a modulus pair $M$, we write
\begin{equation}\label{eq:ZtrMod}
 \Ztr(M):=\MCor(-,M)\in\MPST
\end{equation}
for the represented modulus presheaf; see \cite[Section~3.5]{RSY22}.
See also \cite[Definition~5.1]{HS25} for the conventions on modulus pairs.

Throughout the paper, $\otimes$ denotes the tensor product induced by the
symmetric monoidal structure of \cite[Equation~(2.5.1)]{RSY22}.  On
modulus pairs it is given by
\begin{equation}\label{eq:sum}
 (\overline M,M^\infty)\otimes(\overline N,N^\infty)
 =
 (\overline M\times\overline N,
 M^\infty\times\overline N+\overline M\times N^\infty),
\end{equation}
which is \cite[Equation~(2.5.1)]{RSY22}.  The cube is
$\bcube=(\bP^1,\infty)$.  
A presheaf $\cF\in \MPST$ is said to be $\bcube$-invariant if 
\[
\mathrm{pr}_1^*: \cF(M)\to\cF(M\otimes \bcube)
\]
is an isomorphism for all $M\in \MCor$.
Denote by $\CI$ the full subcategory of $\bcube$-invariant presheaves.
For $\mathcal F\in\MPST$, its maximal
$\bcube$-invariant quotient is
\begin{equation}\label{eq:h0c}
 h_0^{\bcube}(\mathcal F)(M)
 :=\Coker\left(
 \mathcal F(M\otimes\bcube)
 \xrightarrow{i_0^*-i_1^*}
 \mathcal F(M)
 \right).
\end{equation}
This is the maximal $\bcube$-invariant quotient; compare
\cite[Definition~5.3 and Remark~5.4]{HS25}.  For a proper modulus pair
$M$, we write
\begin{equation}\label{eq:h0pair}
 h_0^{\bcube}(M):=h_0^{\bcube}(\Ztr(M)),
\end{equation}
as in \cite[Section~3.6, especially Equation~(3.6.1)]{RSY22}.  In
particular, $h_0^{\bcube}(\Spec k,\varnothing)=\Z$.
For a proper modulus pair $M$, put
\begin{equation}\label{eq:h0red}
 h_0^{\bcube}(M)^0
 :=\Ker\left(
 h_0^{\bcube}(M)\to h_0^{\bcube}(\Spec k,\varnothing)=\Z
 \right),
\end{equation}
where the map is induced by the structure morphism of $M$, as in
\cite[Equation~(5.1.1)]{RSY22}.

Let $\omega:\MCor\to\Cork$ be the interior functor
$M\mapsto M^\circ$.  Precomposition with $\omega$ gives a functor
$\omega^*:\PST\to\MPST$, explicitly
\begin{equation}\label{eq:omstar}
 (\omega^*\mathcal G)(M)=\mathcal G(M^\circ).
\end{equation}
We denote by $\omega_!:\MPST\to\PST$ its left adjoint; equivalently,
$\omega_!$ is the left Kan extension along $\omega$.  Thus
\begin{equation}\label{eq:ombang}
 \Hom_{\PST}(\omega_!\mathcal F,\mathcal G)
 \simeq
 \Hom_{\MPST}(\mathcal F,\omega^*\mathcal G).
\end{equation}
See \cite[Sections~1.3 and~3.5]{RSY22}.  We write
\begin{equation}\label{eq:h0nis}
 h_0(\mathcal F):=\omega_!h_0^{\bcube}(\mathcal F),
 \qquad
 h_{0,\Nis}(\mathcal F):=a_{\Nis}h_0(\mathcal F),
\end{equation}
where $a_{\Nis}$ denotes Nisnevich sheafification.  For a proper modulus
pair $M$, we also write $h_0(M):=h_0(\Ztr(M))$.  The canonical map
$\Ztr(M)\to h_0^{\bcube}(M)$ induces, after applying $\omega_!$, a
canonical surjection
\begin{equation}\label{eq:Ztrh0}
 \Ztr(M^\circ)\simeq\omega_!\Ztr(M)\to h_0(M).
\end{equation}
These are the conventions used in
\cite[Section~3.10 and Equation~(3.15.1)]{RSY22}.

We write $\RSC$ for the full subcategory of $\PST$ of reciprocity
presheaves and $\RSC_{\Nis}$ for its Nisnevich sheaves; see
\cite[Definition~5.7]{HS25} and \cite[Section~3.10]{RSY22}.
The functor $\omega_!$ restricts to the exact localization
\begin{equation}\label{eq:omCI}
 \omega_{\CI}:\CI\to\RSC
\end{equation}
of \cite[Section~3.11, Equation~(3.11.1)]{RSY22}.

\subsection{Cube-invariant objects for the additive and the multiplicative group}
Put
$\Gapair :=(\bP^1,2\infty),\ \Gmpair:=(\bP^1,0+\infty)\in \MCor$.
Let $p:\Gapair \to(\Spec k,\varnothing), q:\Gmpair \to(\Spec k,\varnothing)$ be the structure morphisms, and let
$i_0:(\Spec k,\varnothing)\to\Gapair$ (resp.~$i_1:(\Spec k,\varnothing)\to\Gmpair$) be the zero section (resp.~the unit section).  
We use the pointed modulus presheaves;
\begin{equation}\label{eq:GaM0}
 \GaMod
 :=\Coker\left(
 \Z=\Ztr(\Spec k,\varnothing)
 \xrightarrow{\,\Ztr(i_0)\,}
 \Ztr(\Gapair)
 \right)\in\MPST, 
\end{equation}
and 
\begin{equation}\label{eq:GmM0}
 \GmMod
 :=\Coker\left(
 \Z=\Ztr(\Spec k,\varnothing)
 \xrightarrow{\,\Ztr(i_1)\,}
 \Ztr(\Gmpair)
 \right)\in\MPST.
\end{equation}
Note that the objects $\GaMod, \GmMod$ are denoted by $\GaMod/0, \GmMod/1$ respectively in \cite[Proposition~5.6 and Theorem~5.14]{RSY22}.
Since
$p\circ i_0=\Id$ and $q\circ i_1=\Id$, we have canonical splittings
\begin{align}\label{eq:splitGaGmMod}
    \Ztr(\Gapair)\simeq \Z\oplus \GaMod,\quad 
    \Ztr(\Gmpair)\simeq \Z\oplus \GmMod.
\end{align}
Let
\begin{equation}\label{eq:GaObj}
  \GaSharp:=h_0^{\bcube}(\Gapair)^0,\quad \GmSharp:=h_0^{\bcube}(\Gmpair)^0\in\CI.
\end{equation}
Thus $\GaSharp, \GmSharp$ are exactly the object denoted by
$\mathbb G_a^{\#}, \GmSharp$ in \cite[Section~5.1]{RSY22}, and
\cite[Equation~(5.1.2)]{RSY22} gives the explicit identification
$$ 
a_{\Nis}\omega_{\CI}(\GaSharp)\simeq\Ga,\quad a_{\Nis}\omega_{\CI}(\GmSharp)\simeq\Gm.
$$
By \cite[Section~5.1]{RSY22},
$$h_0^{\bcube}(\GaMod)\simeq\GaSharp,\quad
h_0^{\bcube}(\GmMod)\simeq\GmSharp.
$$
Applying $h_0^{\bcube}$ gives the canonical splittings;
\begin{equation}\label{eq:splitGaGm}
 h_0^{\bcube}(\Gapair)\simeq \Z\oplus\GaSharp,
 \quad 
 h_0^{\bcube}(\Gmpair)\simeq \Z\oplus\GmSharp.
\end{equation}
By \cite[Remark~3.8]{RSY22}, $\CI$ is symmetric monoidal for this
tensor product and $h_0^{\bcube}$ is monoidal.  Hence, for $r\geq1$, fix $*\in\{a,m\}$,
\begin{equation}\label{eq:monGa}
 h_0^{\bcube}((\mathcal{G}_*)^{\otimes r})
 \simeq
 h_0^{\bcube}(\mathcal{G}_*)^{\overset{\CI}{\otimes} r}
 \simeq
 (\Z\oplus \mathbb{G}_*^{\#})^{\overset{\CI}{\otimes} r}.
\end{equation}
In particular, $(\mathbb{G}_*^{\#})^{\overset{\CI}{\otimes}r}$ is a canonical direct
summand of $h_0^{\bcube}((\mathcal{G}_*)^{\otimes r})$.

For $r\geq1$ and $s\geq0$, put
\begin{equation}\label{eq:Ars}
 \Acal_{r,s}:=
 h_{0,\Nis}\left(
 (\GaSharp)^{\otimes r}\otimes(\GmSharp)^{\otimes s}
 \right)\in \RSC_{\Nis}.
\end{equation}
By the monoidality above, $\Acal_{r,s}$ is canonically identified with
$a_{\Nis}\omega_{\CI}((\GaSharp)^{\otimes r}\otimes(\GmSharp)^{\otimes s})$.  Thus the sum $r+s$ of the subscripts
of $\Acal_{r,s}$ is always the number of tensor factors.

Iterating
\cite[Lemma~3.9]{RSY22} therefore gives, for every $r\geq1$ and $s\ge 0$, a canonical
isomorphism
\begin{equation}\label{eq:GaM0Ar}
 h_{0,\Nis}\bigl((\GaMod)^{\otimes r}\otimes (\GmMod)^{\otimes s}\bigr)
 \simeq \Acal_{r,s}.
\end{equation}
This is the comparison used in the direct construction of the differential
morphisms $\Acal_{r,s}\to \Omega^{q+s}\ (0\le q\le r-1)$ below.

\subsection{Presentation of the tensor products}
Let $F/k$ be a field. For $r\geq1$ and $s\geq  0$, put
\begin{equation}\label{eq:symbol-group}
 \mathcal S_{r,s}(F)
 :=\bigoplus_{E/F\,\mathrm{finite}}
 E^{\otimes_{\Z}r}\otimes_{\Z}(E^\times)^{\otimes_{\Z}s},
\end{equation}
where each $E$ is regarded as its additive group.
For a finite extension
$E/F$ and  $a_i\in E$ and $b_j\in E^\times$, the element corresponding to $a_1\otimes \cdots \otimes a_r\otimes b_1\otimes\cdots\otimes b_s$ in the $E$-summand of \eqref{eq:symbol-group} is written
\begin{equation}\label{eq:field-symbol}
 [a_1,\ldots,a_r;b_1,\ldots,b_s]_{E/F}=[\va;\vb]_{E/F}.
\end{equation}
Here $\va$ (resp.~$\vb$) denote the $a_i$-parts (resp.~$b_i$-parts).
If $E=F$, we abbreviate this
to $[a_1,\ldots,a_r;b_1,\ldots,b_s]_F$.
Thus the bracket is $\Z$-additive in every
entry; for example in $s=0$,
\begin{equation}\label{eq:addsym}
 [a_1+a_1',a_2,\ldots,a_r]_{E/F}
 = [a_1,a_2,\ldots,a_r]_{E/F}
 +[a_1',a_2,\ldots,a_r]_{E/F},
\end{equation}
and the analogous equality holds in every other entry.  Since $\operatorname{char}(k)=p$, 
if $p>0$, then we also have
\begin{equation}\label{eq:ptor}
 p[a_1,\ldots,a_r;b_1,\ldots,b_s]_{E/F}
 =[pa_1,a_2,\ldots,a_r;b_1,\ldots,b_s]_{E/F}=0.
\end{equation}
Thus every quotient of $\mathcal S_{r,s}(F)$ considered below is killed by $p$.

Let $F\subset E\subset L$ be finite extensions.
The
\emph{projection-formula relation} is the following;
\begin{enumerate}
    \item Fix $1\leq j\leq r$,
    let $c\in L$, and let $a_i\in E$ for $i\neq j$, let $\vb\in (E^\times)^s$.
    \begin{equation}\label{eq:PF}
    [a_1,\ldots,a_{j-1},\Tr_{L/E}(c),a_{j+1},\ldots,a_r;\vb]_{E/F} =
    [a_1,\ldots,a_{j-1},c,a_{j+1},\ldots,a_r;\vb]_{L/F},
    \end{equation}
    where the entries $a_i$ on the right are viewed in $L$.  
    
    \item 
    Fix $1\leq j\leq s$,
    let $c\in L^\times$, and let $b_i\in E^\times$ for $i\neq j$, let $\va\in E^r$.  
    \begin{equation}\label{eq:PF2}
    [\va;b_1,\ldots,b_{j-1},N_{L/E}(c),b_{j+1},\ldots,b_s]_{E/F} =
    [\va;b_1,\ldots,b_{j-1},c,b_{j+1},\ldots,b_s]_{L/F},
    \end{equation}
    where the entries $b_i$ on the right are viewed in $L$. 
\end{enumerate}
We denote by $\textnormal{(PF)}$ the subgroup of $\mathcal S_{r,s}(F)$
generated by the differences of the two sides of \eqref{eq:PF} and \eqref{eq:PF2}.  


Let
\[
 \tau:\MCor\hookrightarrow \ulMCor
\]
be the inclusion.  Following \cite[Equation~(3.5.1)]{RSY22}, we denote by
\[
 \tau_!:\MPST\to\ulMPST
\]
the exact extension functor from presheaves on proper modulus pairs to
presheaves on all modulus pairs.

\begin{dfn}[Weil reciprocity relation]
\label{dfn:WR}
Let $F/k$ be a field and let $C/F$ be a regular projective
curve.  
Let $r\geq1$ and $s\geq 0$. Choose effective divisors $D_1,\ldots,D_{r+s}$ on $C$ and sections
\begin{equation}\label{eq:WRsec}
 a_i\in(\tau_!\GaSharp)(C,D_i)\quad(1\leq i\leq r),
 \qquad 
 b_j\in(\tau_!\GmSharp)(C,D_{r+j})\quad(1\leq j\leq s)
\end{equation}
For the tensor product \eqref{eq:sum}, put
$D:=D_1+\cdots+D_{r+s}$.
Let $f\in F(C)^\times$ satisfy $f\equiv1\pmod D$.  If
$D=\sum_x n_xx$, this means
$v_x(f-1)\geq n_x$ 
for every $x\in|D|$,
where $v_x$ is the normalized valuation at $x$.  For a closed point
$x\in C\smallsetminus|D|$, let
\begin{equation}\label{eq:WRspec}
 a_i(x)\in\GaSharp(F(x))=F(x),\qquad
 b_j(x)\in\GmSharp(F(x))=F(x)^\times
\end{equation}
be the specializations of $a_i$ and $b_j$ at $x$, in the notation of
\cite[Section~4.8]{RSY22}.  The \emph{Weil reciprocity relation}
$\textnormal{(WR)}$ is
\begin{equation}\label{eq:WR}
 \sum_{x\in C\smallsetminus|D|}
 v_x(f)\,
 [a_1(x),\ldots,a_r(x);b_1(x),\dots,b_s(x)]_{F(x)/F}=0
\end{equation}
in $\mathcal S_{r,s}(F)/\textnormal{(PF)}$.  This is precisely the relation
$\textnormal{(R2)}$ of \cite[Theorem~4.9]{RSY22}, specialized to
$F_1=\cdots=F_r=\GaSharp$ and $F_{r+1}=\cdots =F_{r+s}=\GmSharp$.
\end{dfn}

\begin{thm}[{\cite[Theorem~4.9]{RSY22}}]\label{thm:pres}
Let $F/k$ be a field. Let $r\geq1$ and $s\geq 0$.
The canonical map from the symbol
group \eqref{eq:symbol-group} induces an isomorphism
\begin{equation}\label{eq:pres}
 \frac{\mathcal S_{r,s}(F)}{\textnormal{(PF)},\textnormal{(WR)}}
 \xrightarrow{\;\simeq\;}
 \Acal_{r,s}(F),
\end{equation}
where $\textnormal{(PF)}$ is the projection-formula relation
\eqref{eq:PF} and \eqref{eq:PF2},  
and $\textnormal{(WR)}$ is \autoref{dfn:WR}.
Under this isomorphism the
class of $a_1\otimes\cdots\otimes a_r\otimes b_1\otimes \cdots \otimes b_s$ is denoted by the  symbol
$[a_1,\ldots,a_r;b_1,\ldots,b_s]_{E/F}$.
\end{thm}

\begin{lem}\label{lem:Fvec}
Let $F/k$ be a field and let $r\geq1$ and $s\ge 0$.  The group $\Acal_{r,s}(F)$ has a natural $F$-vector space structure.  
The scalar multiplication is given by
\begin{equation}\label{eq:Fact}
 \lambda [a_1,a_2,\ldots,a_r;\vb]_{E/F}
 :=[\lambda a_1,a_2,\ldots,a_r;\vb]_{E/F}
 \qquad(\lambda\in F).
\end{equation}
\end{lem}
\begin{proof}
For $\lambda\in F$, define an endomorphism
\[
 T_\lambda:\mathcal S_{r,s}(F)\to\mathcal S_{r,s}(F); [a_1,a_2,\ldots,a_r;\vb]_{E/F}\mapsto  [\lambda a_1,a_2,\ldots,a_r;\vb]_{E/F}.
\]
Since the bracket is $\Z$-additive in the first entry, $T_\lambda$ is a
homomorphism of abelian groups. 
Moreover, for $\lambda,\mu\in F$,
\begin{align*}
 T_{\lambda+\mu}
 ([a_1,\ldots,a_r;\vb]_{E/F})
&=(T_\lambda+T_\mu)
 ([a_1,\ldots,a_r;\vb]_{E/F}),\quad   \mbox{and}\\
  T_{\lambda\mu}
 ([a_1,\ldots,a_r;\vb]_{E/F})
 &= T_\lambda T_\mu
 ([a_1,\ldots,a_r;\vb]_{E/F}).
\end{align*} 
Also $T_1=\Id$ and $T_0=0$.  Hence the maps $T_\lambda$ define an
$F$-vector space structure on $\mathcal S_{r,s}(F)$. 
It remains to show that the subgroup generated by $\textnormal{(PF)}$ and $\textnormal{(WR)}$ 
is stable under every $T_\lambda$. 

We first consider \textnormal{(PF)}.  
Let $F\subset E\subset L$ be finite field extensions.
Let $1\leq j\leq r$,
    let $c\in L$, let $a_i\in E$ for $i\neq j$, and let $\vb\in (E^\times)^s$.
If $j=1$, we have 
\begin{align*}
 &[\lambda\Tr_{L/E}(c),a_2,\ldots,a_r;\vb]_{E/F}
 -
 [\lambda c,a_2,\ldots,a_r;\vb]_{L/F}\\
 &\qquad=
 [\Tr_{L/E}(\lambda c),a_2,\ldots,a_r;\vb]_{E/F}
 -
 [\lambda c,a_2,\ldots,a_r;\vb]_{L/F},
\end{align*}
because $\lambda\in F\subset E$ and therefore $\Tr_{L/E}(\lambda c)=\lambda\Tr_{L/E}(c)$.
The last expression is again a projection-formula relation.
Suppose next that $j\neq1$.  Applying $T_\lambda$ to the relation
\[
 [a_1,\ldots,\Tr_{L/E}(c),\ldots,a_r;\vb]_{E/F}
 -
 [a_1,\ldots,c,\ldots,a_r;\vb]_{L/F}
\]
gives
\[
 [\lambda a_1,\ldots,\Tr_{L/E}(c),\ldots,a_r;\vb]_{E/F}
 -
 [\lambda a_1,\ldots,c,\ldots,a_r;\vb]_{L/F}.
\]
This is again a projection-formula relation, now with $\lambda a_1$ in
the first fixed entry.  
Let $1\leq j\leq s$,
    let $c\in L^\times$, let $b_i\in E^\times$ for $i\neq j$, and let $\va\in E^r$.  
    Applying $T_\lambda$ to the relation
\[
 [\va;b_1,\ldots,N_{L/E}(c),\ldots,b_s]_{E/F}
 -
 [\va;b_1,\ldots,c,\ldots,b_s]_{L/F}
\]
gives
\[
 [\lambda a_1,a_2,\ldots,a_r;b_1,\ldots,N_{L/E}(c),\ldots,b_s]_{E/F}
 -
 [\lambda a_1,a_2,\ldots,a_r;b_1,\ldots,c,\ldots,b_s]_{L/F}.
\]
This is again a projection-formula relation.
Thus the subgroup generated by
\textnormal{(PF)} is stable under $T_\lambda$.

We next consider \textnormal{(WR)}.  
Let $C/F$, $D_i$, $a_i$, $b_j$ and $f$
be data as in \autoref{dfn:WR}; thus
$a_i\in(\tau_!\GaSharp)(C,D_i)\ (1\le i\le r)$, $b_j\in(\tau_!\GmSharp)(C,D_{r+j})\ (1\le j\le s)$, 
$D=D_1+\cdots+D_{r+s}$, $f\equiv1\pmod D$, 
and the corresponding relation is
\[
R := \sum_{x\in C\smallsetminus|D|}
 v_x(f)
 [a_1(x),\ldots,a_r(x);b_1(x),\dots, b_s(x)]_{F(x)/F}.
\]
Suppose first that $\lambda\neq0$.  The morphism
\[
 m_\lambda:
 (\mathbb P^1_F,2\infty)
 \to
 (\mathbb P^1_F,2\infty),
 \qquad
 t\mapsto \lambda t,
\]
is an automorphism of modulus pairs.  
Indeed, $m_\lambda^*(2\infty)=2\infty$.
Put $\mathbb G_{a,F}^\#
 =
 h_0^{\overline{\square}}
 (\mathbb P^1_F,2\infty)^0$. 
The map $m_\lambda$ induces an automorphism
\[
 m_\lambda^\#:
 \mathbb G_{a,F}^\#
 \to 
 \mathbb G_{a,F}^\#.
\]
Applying $\tau_{F,!}$ gives an endomorphism
\[
 \tau_{F,!}(m_\lambda^\#):
 \tau_{F,!}\mathbb G_{a,F}^\#
 \to
 \tau_{F,!}\mathbb G_{a,F}^\#.
\]
We therefore define
\[
 \lambda a_1
 :=
 \bigl(\tau_{F,!}(m_\lambda^\#)\bigr)(a_1).
\]
Then
\[
 \lambda a_1
 \in
 \bigl(\tau_{F,!}\mathbb G_{a,F}^\#\bigr)(C,D_1)
 =
 (\tau_!\GaSharp)(C,D_1)\simeq \GaSharp(C,D_1).
\]
For $\lambda=0$, we put $\lambda a_1=0$ the zero of the abelian group $(\tau_!\GaSharp)(C,D_1)$.
For a closed point 
$x\in C\smallsetminus|D_1|$, the specialization map
is functorial with respect to morphisms of modulus presheaves. 
Hence
\[
 (\lambda a_1)(x)
 =
 (m_\lambda^\#)_{F(x)}\bigl(a_1(x)\bigr) =\lambda a_1(x).
\]
Consequently, the \textnormal{(WR)} relation associated with
$\lambda a_1,a_2,\ldots,a_r$ is
\begin{align*}
 &\sum_{x\in C\smallsetminus|D|}
 v_x(f)
 [(\lambda a_1)(x),a_2(x),\ldots,a_r(x);b_1(x),\dots, b_s(x)]_{F(x)/F}\\
 &\qquad=
 \sum_{x\in C\smallsetminus|D|}
 v_x(f)
 [\lambda a_1(x),a_2(x),\ldots,a_r(x);b_1(x),\dots, b_s(x)]_{F(x)/F}\\
 &\qquad=
 T_\lambda
 \left(
 \sum_{x\in C\smallsetminus|D|}
 v_x(f)
 [a_1(x),a_2(x),\ldots,a_r(x);b_1(x),\dots, b_s(x)]_{F(x)/F}
 \right)\\
 &\qquad=T_\lambda(R).
\end{align*}

Hence $T_\lambda$ sends every relation of type \textnormal{(WR)} to
 another relation of type \textnormal{(WR)}.
\color{black}
\end{proof}

\begin{thm}[{\cite[Theorem~5.14 and Theorem~5.15]{RSY22}}]
\label{thm:two}
Assume that $\operatorname{char}(k)\neq2$.  There is an isomorphism of Nisnevich reciprocity
sheaves
\begin{equation}\label{eq:twoiso}
 \Theta_{2,0}:\Acal_{2,0} 
 \xrightarrow{\simeq}
 \Ga\oplus\Omega^1.
\end{equation}
On a field $F/k$, the map is given by
\begin{equation}\label{eq:twofld}
 [x,y]_{E/F}
 \mapsto
 \left(\Tr_{E/F}(xy),\Tr_{E/F}(x\, \d y)\right).
\end{equation}
\end{thm}
In 
\cite[Theorems~5.14 and~5.15]{RSY22}, 
the assertion is obtained, together with the
splitting of the sheaf $\Pcal^1$
if $k$ has characteristic zero. 
However, 
the above isomorphism is valid for characteristic $p\neq 2$.

Write $\sigma_{1,0}:\Omega^1\to\Acal_{2,0}$ for the restriction of $\Theta_{2,0}^{-1}$ to the
second summand of $\Ga\oplus\Omega^1$.  On a field $F/k$ it is given by
\begin{equation}\label{eq:s1}
 \sigma_{1,0,F}(A\, \dlog b)=[A/b,b]_F-[A,1]_F,
 \qquad
 \sigma_{1,0,F}(a\, \d b)=[a,b]_F-[ab,1]_F.
\end{equation}
Since $\sigma_{1,0}$ is a morphism in $\RSC_{\Nis}$, these maps commute with finite
transfers.  

We now sharpen \cite[Theorem~5.19]{RSY22}. 
\begin{thm}\label{thm:one-additive}
If $k$ is a perfect field of $\operatorname{char}(k)\neq 2$, then 
there is a canonical isomorphism in $\RSC_{\Nis}$
\begin{equation}\label{eq:one-additive}
 \Theta_{1,m}:\Acal_{1,m}\xrightarrow{\simeq }\Omega^m
\end{equation}
for $m\geq0$. 
For a finitely generated field $F/k$, a finite extension $E/F$, and
$a\in E$, $b_1,\ldots,b_m\in E^\times$, the field value is
\begin{equation}\label{eq:one-additive-field}
 [a;b_1,\ldots,b_m]_{E/F}
 \mapsto
 \Tr_{E/F}\left(a\,\dlog b_1\wedge\cdots\wedge \dlog b_m\right).
\end{equation}
\end{thm}
\begin{proof}
The case $m=0$ is $\Acal_{1,0}\simeq\Ga$. Hence we assume $m\geq1$.
For $\operatorname{char}(k)\neq3,5$, this is \cite[Theorem~5.19]{RSY22}.  We assume
$\operatorname{char}(k)=p=3$ or $5$.
The morphism
\[
 \Theta_{1,m}:\Acal_{1,m}\to\Omega^m
\]
is constructed in \cite[proof of Theorem~5.19, (5.19.1)]{RSY22}.
The construction is valid also in characteristics $3$ and $5$.
We show that the inverse on fields 
\[ 
\sigma_{0,m,F}:\Omega_F^m\to\Acal_{1,m}(F),
 \qquad
 c\,\dlog b_1\wedge\cdots\wedge \dlog b_m
 \mapsto[c;b_1,\ldots,b_m]_F.
\]
is well-defined and surjective (cf.~\cite[Proposition~5.18]{RSY22} in the case of $\operatorname{char}(F) \neq 3,5$).
For a field $F/k$, consider
\[
 \widetilde \sigma_{0,m,F}:
 F\otimes_{\Z}(F^\times)^{\otimes m}\to\Acal_{1,m}(F),
 \qquad
 c\otimes b_1\otimes\cdots\otimes b_m
 \mapsto[c;b_1,\ldots,b_m]_F.
\]
Since $\Acal_{1,m}(F)$ is an $F$-vector space through the first additive
entry, $\widetilde \sigma_{0,m,F}$ is $F$-linear. 
By \autoref{cor:Aflin}, 
it is enough to check 
\begin{equation}\label{eq:RSY-rep-relation}
[a;b_1,\ldots b_m]_F = 0\quad  \mbox{if $b_i = b_j$ for some $i<j$}, 
\end{equation}
and
\begin{equation}\label{eq:RSY-add-relation}
 [a;a,\vb]_F+[1-a;1-a,\vb]_F=0,
\end{equation}
where $a,1-a\in F^\times$ and $\vb = (b_2,\ldots,b_m)$. 
The first relation \eqref{eq:RSY-rep-relation} is proved in
\cite[proof of Proposition~5.18, part~(ii)]{RSY22}; the argument only
uses that $2$ is invertible.  We prove \eqref{eq:RSY-add-relation}.

Put $A=a(1-a), C=a^2-a+1=1-A$.
First suppose $C=0$. 
If $p=3$, then $C = (a+1)^2$ hence 
 $a=-1$. We have
\[
 [a;a,\vb]_F+[1-a;1-a,\vb]_F
 =
 2[-1;-1,\vb]_F
 =
 [-1;1,\vb]_F=0.
\]
If $p=5$, both $a$ and $1-a$ are roots of $u^2 - u+1$. 
They satisfy $u^6=1$.  
Since $\Acal_{1,m}(F)$ is killed by $5$, for $u=a,1-a$ we have
\[
 [u;u,\vb]_F
 =
 6[u;u,\vb]_F
 =
 [u;u^6,\vb]_F
 =
 [u;1,\vb]_F=0.
\]
Thus \eqref{eq:RSY-add-relation} holds when $C=0$.
Assume $C\neq0$.  On $\bP^1 = \bP_F^1$ with coordinate $t$, take
\begin{equation}\label{eq:new-RSY-function}
 f=
 \frac{(t-a)(t-(1-a))(t^2+1)}
 {(t^2-t+1)(t^2+A)},
 \ 
 g_0=t,\ g_1=t,\ g_i = b_i\ (2\leq i\leq m).
\end{equation}
Since
\[
 (t-a)(t-(1-a))(t^2+1) 
 -(t^2-t+1)(t^2+A) =
 -Ct,
\]
we have
\[
 f-1= \frac{-Ct}{(t^2-t+1)(t^2+A)}.
\]
As $A\neq0$, $v_0(f-1)\geq1$, $v_\infty(f-1)\geq3$.
Hence
$f\equiv1\pmod{(0)+3(\infty)}$. 

Put 
$D_0 = 2(\infty), D_1 = (0)+ (\infty), D_i = 0\ (2\leq i\leq m)$ 
and 
\[
D  = D_0 + D_1 + \cdots + D_m = (0) + 3(\infty).
\]
For $g_0 = t \in F(\bP^1) = F(t)$, 
the morphism  $\Gamma_{g_0}\colon (\bP^1, D_0)\to \Gapair =(\bP^1,2(\infty))$ 
is admissible; 
$D_0 \geq 2\mathrm{div}_\infty(g_0)$. 
We have $g_0 \in h_0^{\bcube}(\Gapair)(\bP^1,D_0)$. 
By the projection $h_0^{\bcube}(\Gapair)\simeq \Z \oplus \GaSharp\to \GaSharp$, we obtain 
\[
g_0^\#\in \GaSharp(\bP^1,2(\infty)) = (\tau_!\GaSharp)(\bP^1,2(\infty)).
\]
In the same way, we have 
\[
 g_1^\#\in(\tau_!\GmSharp)(\bP^1,(0)+(\infty)),\quad g_i^\# \in (\tau_!\GmSharp)(\bP^1,\varnothing)\ (2\le i\le m).
\]
Applying the relation $\mathrm{(WR)}$ of \autoref{thm:pres}, we obtain
\begin{align*}
	 0 &= \sum_{x \in \bP^1\smallsetminus \set{0,\infty}}v_x(f) \Tr_{F(x)/F}([g_0(x);g_1(x), \vb]_{F(x)}) \\
 &=[a;a,\vb]_F+[1-a;1-a,\vb]_F
 +R_1-H-R_A.
\end{align*}
Here, 
for $\lambda = 1$ or $A$, put 
\[
R_\lambda = \sum_{x \mid t^2+ \lambda}v_x(t^2+\lambda)\Tr_{F(x)/F}([g_0(x);g_1(x), \vb]_{F(x)})
\]
where $x$ runs over the closed points of $\bP^1\smallsetminus \set{0,\infty}$ 
corresponding to the irreducible factors of $t^2+\lambda$. 
Since $t^2 + \lambda$ is separable, $v_x(t^2+ \lambda) = 1$. 
In the same way, 
\[
H = \sum_{ x\mid (t^2 - t+1)}v_x(t^2-t+1) \Tr_{F(x)/F}([g_0(x);g_1(x), \vb]_{F(x)}).
\]
We claim that $R_\lambda=0$ for $\lambda =1$ or $A$. 
For each closed point $x$ of $\bP^1 \smallsetminus \set{0,\infty}$ 
corresponding to the irreducible factors of $t^2+\lambda$, 
let $\alpha_x \in F(x)$ be the image of $t$ in the residue field $F(x)$.  
Then $\alpha_x^2 = -\lambda$. 
Thus, by the projection formula,
\[
 2R_\lambda=
 \sum_{x\mid t^2+ \lambda} 
 [\alpha_x;\alpha_x^2,\vb]_{F(x)/F} 
 =
 \left[
 \sum_{x\mid t^2 + \lambda}\Tr_{F(x)/F}(\alpha_x);
 -\lambda,\vb
 \right]_F.
\]
Since the coefficient of $t$ in $t^2+\lambda$ is zero, 
we have $\sum_{x\mid t^2+\lambda} \Tr_{F(x)/F}(\alpha_x) = 0.$
$\Acal_{1,m}(F)$ is killed by the odd prime $p$, multiplication by $2$
is injective, and therefore $R_\lambda=0$.
It remains to treat $H$. 
If $p=3$, then $t^2-t+1=(t+1)^2$,
so
\[
 H=2[-1;-1,\vb]_F=[-1;1,\vb]_F=0.
\]
If $p=5$, every root $\gamma$ of $t^2-t+1$ satisfies
$\gamma^6=1$.  Since the relevant symbol group is killed by $5$,
\[
 [\gamma;\gamma,\vb]
 =
 6[\gamma;\gamma,\vb]
 =
 [\gamma;\gamma^6,\vb]
 =
 [\gamma;1,\vb]
 =0.
\]
Hence $H=0$ also in characteristic $5$.
Therefore the equality 
\eqref{eq:RSY-add-relation} holds.  
Thus $\widetilde  \sigma_{0,m,F}$ induces
\[
 \sigma_{0,m,F}:\Omega_F^m\to\Acal_{1,m}(F),
 \qquad
 c\,\dlog b_1\wedge\cdots\wedge \dlog b_m
 \mapsto[c;b_1,\ldots,b_m]_F.
\]
The compatibility with finite traces is proved as in
\cite[proof of Proposition~5.18, (5.18.3)--(5.18.4)]{RSY22} 
and is valid for every $p\neq 2$.
The rest of the proofs of
\cite[Proposition~5.18 and Theorem~5.19]{RSY22} applies unchanged.
Finally, \cite[Lemma~4.16]{RSY22} gives the required isomorphism in
$\RSC_{\Nis}$.
\end{proof}

\subsection{Auxiliary lemmas}

\begin{lem}\label{lem:fixed-insertion}
Let $F/k$ be a field.  Let $r,r'\geq1$ and $s,s'\geq0$ with
$r\leq r'$ and $s\leq s'$.  Choose 
$1\leq i_1<\cdots<i_r\leq r'$
and 
$1\leq j_1<\cdots<j_s\leq s'$.
Fix $\lambda_i\in F$ for $i\notin\set{i_1,\ldots,i_r}$ and
$\mu_j\in F^\times$ for $j\notin\set{j_1,\ldots,j_s}$.
Then the assignment 
$\mathcal S_{r,s}(F) \to \mathcal S_{r',s'}(F)$
on symbols
\begin{equation}\label{eq:fixed-insertion}
 [a_1,\ldots,a_r;b_1,\ldots,b_s]_{E/F}
 \mapsto
 [c_1,\ldots,c_{r'};d_1,\ldots,d_{s'}]_{E/F},
\end{equation}
where
\[
 c_{i_\nu}=a_\nu,\qquad
 c_i=\lambda_i\quad(i\notin\set{i_1,\ldots,i_r}),
\]
and
\[
 d_{j_\nu}=b_\nu,\qquad
 d_j=\mu_j\quad(j\notin\set{j_1,\ldots,j_s}),
\]
induces a homomorphism
\[
 \Acal_{r,s}(F)\to\Acal_{r',s'}(F).
\]
\end{lem}

\begin{proof}
For any integer $n\ge 0$, put $[n]:=\set{1,\ldots,n}$ and $[0]:=\varnothing$.
Let
 $I:=\set{i_1,\ldots,i_r}\subset[r']$,  
 $J:=\set{j_1,\ldots,j_s}\subset[s']$.
Let $\widetilde\iota:
 \mathcal S_{r,s}(F)\to\mathcal S_{r',s'}(F)$ 
be the homomorphism defined by \eqref{eq:fixed-insertion}.
We show that $\widetilde\iota$ preserves
\textnormal{(PF)} and \textnormal{(WR)}.

We first consider \textnormal{(PF)}.
Let $F\subset E\subset L$ be finite field extensions.

Suppose first that the projection-formula relation occurs in the
$\nu$-th additive entry of $\mathcal S_{r,s}(F)$.
Thus, for $a_\nu\in L$, $a_\ell\in E$ for $\ell\neq\nu$, and
$b_1,\ldots,b_s\in E^\times$, we have
\[
 [a_1,\ldots,\Tr_{L/E}(a_\nu),\ldots,a_r;
 \vb]_{E/F}
 -
 [a_1,\ldots,a_\nu,\ldots,a_r;
 \vb]_{L/F}, 
\]
where $\vb = (b_1,\ldots,b_s)$. 
Its image under $\widetilde\iota$ is
\[[c_1,\ldots,
   \overset{i_\nu}{\Tr_{L/E}(a_\nu)},
   \ldots,c_{r'};\vd]_{E/F}-
 [c_1,\ldots,
   \overset{i_\nu}{a_\nu},
   \ldots,c_{r'};\vd]_{L/F},
\]
where 
 $c_i=\lambda_i\ (i\notin I),
 d_j=\mu_j\ (j\notin J)$, $\vd = (d_1,\ldots,d_{s'})$,
and the remaining variable entries are the corresponding
$a_\ell$ and $b_\ell$.
Since $\lambda_i\in F\subset E$ and
$\mu_j\in F^\times\subset E^\times$, all fixed entries are defined
over $E$.  Hence the above expression is exactly a projection-formula
relation in the $i_\nu$-th entry of
$\mathcal S_{r',s'}(F)$.

Similarly, suppose that the projection-formula relation occurs in the
$\nu$-th multiplicative entry of $\mathcal S_{r,s}(F)$.  For $b_\nu\in L^\times$, its image is
\[
 [\vc;
 d_1,\ldots,
 \overset{j_\nu}{N_{L/E}(b_\nu)},
 \ldots,d_{s'}]_{E/F} -
 [\vc;
 d_1,\ldots,
 \overset{j_\nu}{b_\nu},
 \ldots,d_{s'}]_{L/F},
\]
where $\vc = (c_1,\ldots,c_{r'})$, 
which is the projection-formula relation in the $j_\nu$-th
 entry of $\mathcal S_{r',s'}(F)$.
Thus $\widetilde\iota$ preserves \textnormal{(PF)}.

We next consider \textnormal{(WR)}.
Let $C/F$ be a regular projective curve and choose effective divisors
$ D_1,\ldots,D_r,D_{r+1},\ldots,D_{r+s}$
and sections $a_\nu\in(\tau_!\GaSharp)(C,D_\nu)$ 
 $(1\leq\nu\leq r)$,
and
 $b_\nu\in(\tau_!\GmSharp)(C,D_{r+\nu})$ ($1\leq\nu\leq s$).
Put
\[
 D:=D_1+\cdots+D_{r+s},
\]
and let $f\in F(C)^\times$ satisfy $f\equiv1\pmod D$.
The corresponding Weil reciprocity relation is
\[
 R=
 \sum_{x\in C\smallsetminus|D|}
 v_x(f)
 [a_1(x),\ldots,a_r(x);
 b_1(x),\ldots,b_s(x)]_{F(x)/F}.
\]

For each $i\notin I$, regard $\lambda_i\in F$ as a constant rational
function on $C$.  Since it has no pole, the constant map
$C\to\bP^1_F$; $x\mapsto\lambda_i$,
is admissible from $(C,\varnothing)$ to
$(\bP^1_F,2\infty)$.
Hence $\lambda_i$ defines a section
$\lambda_i\in(\tau_!\GaSharp)(C,\varnothing)$.
Likewise, for each $j\notin J$, the constant
$\mu_j\in F^\times$ defines a section
$\mu_j\in(\tau_!\GmSharp)(C,\varnothing)$,
because the constant map
$C\to\bP^1_F,\ x\mapsto\mu_j$,
does not meet $0$ or $\infty$.

Define sections $c_1,\ldots,c_{r'}$ and
$d_1,\ldots,d_{s'}$ by
 $c_{i_\nu}=a_\nu$,  $c_i=\lambda_i$ ($i\notin I$),
and
 $d_{j_\nu}=b_\nu$, 
 $d_j=\mu_j$ ($j\notin J$).
The function $f$ satisfies the required congruence condition for these sections.
The corresponding Weil reciprocity relation in
$\mathcal S_{r',s'}(F)$ is
\[
 \sum_{x\in C\smallsetminus|D|}
 v_x(f)
 [c_1(x),\ldots,c_{r'}(x);
 d_1(x),\ldots,d_{s'}(x)]_{F(x)/F}.
\]
Since specialization of a constant section is constant, we have
 $c_{i_\nu}(x)=a_\nu(x)$, 
 $c_i(x)=\lambda_i$ ($i\notin I$),
and
 $d_{j_\nu}(x)=b_\nu(x)$, 
 $d_j(x)=\mu_j$ ($j\notin J$).
Hence this relation is exactly $\widetilde\iota(R)$.
Thus $\widetilde\iota$ preserves \textnormal{(WR)}.

It follows from \autoref{thm:pres} that
$\widetilde\iota$ induces a homomorphism
\[
 \Acal_{r,s}(F)\to\Acal_{r',s'}(F),
\]
whose value on symbols is given by \eqref{eq:fixed-insertion}.
\end{proof}
\black

\begin{lem}\label{lem:tworel}
Let $F/k$ be a field.
\begin{enumerate}
    \item 
    Let $r\geq2$ and $s\ge 0$.  Fix
$1\leq i<j\leq r$.
Then, in $\Acal_{r,s}(F)$, for $x,y,z\in F$, one has
\begin{align}
 &[a_1,\ldots,a_{i-1},x,a_{i+1},\ldots,
   a_{j-1},yz,a_{j+1},\ldots,a_r;\vb]_F
 \notag\\
 &\qquad=
 [a_1,\ldots,a_{i-1},xy,a_{i+1},\ldots,
   a_{j-1},z,a_{j+1},\ldots,a_r;\vb]_F
 \notag\\
 &\qquad\quad+
 [a_1,\ldots,a_{i-1},xz,a_{i+1},\ldots,
   a_{j-1},y,a_{j+1},\ldots,a_r;\vb]_F
 \notag\\
 &\qquad\quad-
 [a_1,\ldots,a_{i-1},xyz,a_{i+1},\ldots,
   a_{j-1},1,a_{j+1},\ldots,a_r;\vb]_F,
 \label{eq:leib}
\end{align}
and
\begin{align}
 &[a_1,\ldots,a_{i-1},x,a_{i+1},\ldots,
   a_{j-1},y,a_{j+1},\ldots,a_r;\vb]_F
 \notag\\
 &\qquad+
 [a_1,\ldots,a_{i-1},y,a_{i+1},\ldots,
   a_{j-1},x,a_{j+1},\ldots,a_r;\vb]_F
 \notag\\
 &\qquad=
 [a_1,\ldots,a_{i-1},1,a_{i+1},\ldots,
   a_{j-1},xy,a_{j+1},\ldots,a_r;\vb]_F
 \notag\\
 &\qquad\quad+
 [a_1,\ldots,a_{i-1},xy,a_{i+1},\ldots,
   a_{j-1},1,a_{j+1},\ldots,a_r;\vb]_F.
 \label{eq:sym}
\end{align}
\item 
Let $r\geq2$ and $s\ge 1$.  Fix
$1\leq i\leq r$ and $1\leq j\leq s$. 
Then, in $\Acal_{r,s}(F)$, for $x,1-x\in F^\times$, one has
\begin{align}
    &[a_1,\ldots,a_{i-1},x,a_{i+1},\ldots,a_r;b_1,\dots,b_{j-1},x,b_{j+1},\dots,b_s]_F 
 \notag\\
 &\ +
 [a_1,\ldots,a_{i-1},1-x,a_{i+1},\ldots,a_r;b_1,\dots,b_{j-1},1-x,b_{j+1},\dots,b_s]_F
 =0.
 \label{eq:diff}
\end{align}
Moreover, for $1\leq j<j'\leq s$,
\begin{align}
 [\va;b_1,\dots,b_{j-1},x,b_{j+1},\dots,b_{j'-1},x,b_{j'+1},\dots,b_s]_F
 =0.
 \label{eq:diff2}
\end{align}
\end{enumerate}
\end{lem}

\begin{proof}
\noindent (1)
First, we show the assertion for $r=2$ and $s=0$. 
By \autoref{thm:two}, one has
\begin{align}\label{eq:leib0}
 [x,yz]_F
 =
 [xy,z]_F+[xz,y]_F-[xyz,1]_F
\end{align}
because under \eqref{eq:twofld} its degree-zero component is the identity
$xyz=xyz+xyz-xyz$, while its degree-one component is the Leibniz rule
$x\d(yz)=xy\,\d z+xz\,\d y$.
Likewise,
\begin{align}\label{eq:sym0}
     [x,y]_F+[y,x]_F=[1,xy]_F+[xy,1]_F,
\end{align}
since the first components in $\Ga(F)$ 
on both sides are $2xy$ and the
second components in $\Omega^1_F$ agree by
$x\,\d y+y\,\d x=\d(xy)$.

Consider the homomorphism
\begin{equation}\label{eq:insert-two}
 \iota_{i,j}:\mathcal S_{2,0}(F)\to\mathcal S_{r,s}(F)
\end{equation}
which sends a symbol $[u,v]_{E/F}$ to
\[
 [a_1,\ldots,a_{i-1},u,a_{i+1},\ldots,
   a_{j-1},v,a_{j+1},\ldots,a_r;\vb]_{E/F}.
\]
By \autoref{lem:fixed-insertion}, 
$\iota_{i,j}$ induces a homomorphism
\begin{equation}\label{eq:insert-two-Ar}
 \iota_{i,j}:\Acal_{2,0}(F)\to\Acal_{r,s}(F).
\end{equation}

Applying $\iota_{i,j}$ to \eqref{eq:leib0} and \eqref{eq:sym0} gives the desired equations \eqref{eq:leib} and \eqref{eq:sym} in $\Acal_{r,s}(F)$.

\smallskip
\noindent (2)
From the proof of \autoref{thm:one-additive}, one has
\begin{align}\label{eq:1}
    [x;x]_F+[1-x;1-x]_F=0 \quad \text{in $\Acal_{1,1}(F)$}
\end{align}
and
\begin{align}\label{eq:2}
    [a;x,x]_F=0 \quad \text{in $\Acal_{1,2}(F)$}.
\end{align}

Consider the homomorphisms
\begin{equation}\label{eq:insert-two2}
 \iota'_{i,j}:\mathcal S_{1,1}(F)\to\mathcal S_{r,s}(F),\quad 
\iota''_{j,j'}: \mathcal{S}_{1,2}(F)\to \mathcal{S}_{r,s}(F)
\end{equation}
such that
\[
 \iota'_{i,j}([u;v]_{E/F})=[a_1,\ldots,a_{i-1},u,a_{i+1},\ldots,a_r;b_1,\dots,b_{j-1},v,b_{j+1},\dots,b_s]_{E/F},
\]
and 
\[
 \iota''_{j,j'}([u;v,w]_{E/F})=[u,a_2,\ldots,a_r;b_1,\dots,b_{j-1},v,b_{j+1},\dots,b_{j'-1},w,b_{j'+1},\dots,b_s]_{E/F}.
\]
By \autoref{lem:fixed-insertion},   
the maps $\iota'_{i,j}$ and $\iota''_{j,j'}$ induce homomorphisms
\begin{equation}\label{eq:insert-mixed-Ar}
 \iota'_{i,j}:\Acal_{1,1}(F)\to\Acal_{r,s}(F),\quad \iota''_{j,j'}: \Acal_{1,2}(F)\to \Acal_{r,s}(F).
\end{equation}
Applying $\iota'_{i,j}$ (resp.~$\iota''_{j,j'}$) to \eqref{eq:1} (resp.~\eqref{eq:2}) gives the desired equation \eqref{eq:diff} (resp.~\eqref{eq:diff2}) in $\Acal_{r,s}(F)$.
\end{proof}

\section{Main Result}
\label{sec:diff}

\subsection{Construction of the differential morphism}
Let $n, m\geq0$.  
Let $t_0,\ldots,t_n$ be the
affine coordinates on the factors of $(\Gapair^{\otimes(n+1)})^\circ = (\bA^1)^{n+1}$.
Let $x_1,\ldots, x_m$ be the
affine coordinates on the factors of $(\Gmpair^{\otimes m})^\circ = (\bA^1\smallsetminus\{0\})^{m}$.
Set
\begin{equation}\label{eq:eta}
 \eta_{n,m}:=
 t_0
 \bigwedge_{i=1}^n\d t_i
 \bigwedge_{j=1}^m\dlog x_j
 \in\Omega^{n+m}((\mathbb A^1)^{n+1}\times (\bA^1\smallsetminus\{0\})^m).
\end{equation}\black
Each $t_i$ and $x_j$ give morphisms
\[
t_i:\Ztr(\bA^1)\to \Ga,\qquad x_j:\Ztr(\bA^1\smallsetminus\{0\})\to \Gm
\]
which descend to 
\[
t_i:\omega_!\Ztr(\Gapair) \to \Ga,\qquad
x_j:\omega_!\Ztr(\Gmpair) \to \Gm.
\]
Furthermore, the composition
\[
\omega_!\Ztr(\Spec k, \varnothing)\xrightarrow{\Ztr(i_0)} \omega_!\Ztr(\Gapair) \xrightarrow{t_i} \Ga
\]
and
\[
\omega_!\Ztr(\Spec k, \varnothing)\xrightarrow{\Ztr(i_1)} \omega_!\Ztr(\Gmpair) \xrightarrow{x_j} \Gm
\]
are trivial, and hence we have morphisms
\[
t_i:\omega_!\GaMod \to \Ga, \qquad 
x_j:\omega_!\GmMod \to \Gm.
\]
Taking tensor product and composing with $\eta_{n,m}$, we have a morphism 
\[
\widetilde{\theta}_{n,m}: \omega_!((\GaMod)^{\otimes (n+1)}\otimes (\GmMod)^{\otimes m}))\xrightarrow{t_0\otimes \cdots \otimes t_n\otimes x_1\otimes\cdots\otimes x_m} \Ga^{\otimes (n+1)}\otimes \Gm^{\otimes m}\xrightarrow{\eta_{n,m}} \Omega^{n+m}.
\]
\begin{lem}\label{lem:thetaS}
    The morphism $\widetilde{\theta}_{n,m}$ factors through $h_0((\GaMod)^{\otimes (n+1)}\otimes (\GmMod)^{\otimes m})$.
\end{lem}

Since $\Omega^{n+m}$ is a Nisnevich sheaf, $\widetilde{\theta}_{n,m}$ factors through $h_{0,\Nis}((\GaMod)^{\otimes (n+1)}\otimes (\GmMod)^{\otimes m})$. From \eqref{eq:GaM0Ar}, we get the desired morphism
\begin{equation}\label{eq:thetaS}
 \theta_{n,m}:\Acal_{n+1,m}\to\Omega^{n+m}.
\end{equation}

\begin{proof}[Proof of \autoref{lem:thetaS} ]
If $n=m=0$, then 
$\eta_{0,0}=t_0$.
Hence $\widetilde{\theta}_{0,0}$ is the canonical morphism
$\omega_!\GaMod\to\Ga$.
By \cite[Section~5.1]{RSY22}, this morphism factors through $h_0(\GaMod)$. 
Thus the assertion holds when $n=m=0$.
We may therefore assume that $n+m>0$.
    
    Via the automorphism $\varphi:\bP^1\to\bP^1; t\mapsto\frac{t}{t-1}$, the modulus pair $\bcube=(\bP^1,\infty)$ is identified with $(\bP^1,1)$, and the sections $0,1$ correspond to $0,\infty$.
    Thus it is enough to prove the following. For any connected $X\in \Sm$ and any prime correspondence $\alpha\in \ulMCor(X\otimes (\bP^1,1), (\bP^1,2\infty)^{\otimes (n+1)}\otimes (\bP^1,0+\infty)^{\otimes m})$,
    \[
    (i_0^\ast-i_\infty^\ast)\alpha^\ast (\eta_{n,m})=0 \ \text{in}\ \Omega^{n+m}(X).
    \]
    Since the sheaf $\Omega^{n+m}$ has global injectivity, the canonical map $\Omega^{n+m}(X)\to \Omega^{n+m}(k(X))$ is injective.
    Therefore we may assume that $X=\Spec F$ where $F$ is a function field.
    Hence the correspondence $\alpha$ can be written as
    \[
    \alpha=\pi\circ \Gamma_f^t
    \]
    where $\pi: C\to (\bP^1)^{n+1}\times (\bP^1)^m$ is a morphism from a proper regular curve over $F$, $f\in F(C)^\times$ satisfying $f\equiv 1 \bmod D_C:=\pi^*(D)$ and $\Gamma_f^t\subset \bP^1\times C$ is the transpose of the graph of $f$.
    Here $D$ is the effective divisor of the modulus pair $(\bP^1,2\infty)^{\otimes (n+1)}\otimes (\bP^1,0+\infty)^{\otimes m}$. 
    Then we have
    \begin{align}\label{eq:res}
        (i_0^\ast-i_\infty^\ast)\alpha^\ast \eta_{n,m}
        &=\mathrm{div}_C(f)^\ast \pi^\ast(\eta_{n,m})\\
        &=\sum_{x\in C\smallsetminus |D_C|}\mathrm{Res}_{C,x}(\pi^\ast(\eta_{n,m})\dlog f ),
    \end{align}
    where $\mathrm{Res}_{C,x}$ is the residue map for absolute K\"ahler differentials. 
    Let $p_i:(\bP^1)^{n+1}\times (\bP^1)^m\to \bP^1$ be the projection to the $i$-th factor.
We write $m_i$ for the multiplicity of $x\in|D_C|$ in $\pi^*p_i^*(\infty)$ and put $m_j':=v_x(\pi^*x_j)$.
Put $M:=\sum_{i=0}^n m_i$.
Let $u$ be a local parameter at $x$.
Then $\pi^*t_i$, $\pi^*x_j$ and $f$ can be written as follows:
\[
\pi^*t_i=a_i u^{-m_i},\quad \pi^*x_j=b_j u^{m_j'},\quad f=1+cu^e,
\]
where $a_i,c\in\mathcal{O}_{C,x}$, $b_j\in\mathcal{O}_{C,x}^\times$ and $e\geq 2M+\sum_{j=1}^m |m_j'|$. 
    Then, we obtain that 
    \begin{align}
        \pi^\ast(\eta_{n,m})
        &=a_0u^{-m_0}\bigwedge_{i=1}^n(u^{-m_i}\d a_i-m_ia_iu^{-m_i-1}\d u)\bigwedge_{j=1}^m(\dlog b_j+m_j'\dlog u)\\
        &=u^{-1-M}\omega_1\wedge \d u+u^{-M}\omega_2
    \end{align}
    where $\omega_1\in \Omega^{n+m-1}_{C,x}$ and $\omega_2\in \Omega_{C,x}^{n+m}$.
    On the other hand, we have
    \[
    \dlog f=\frac{1}{f}(u^e\d c+ecu^{e-1}\d u).
    \]
    Thus, we have 
    \[
    \pi^\ast(\eta_{n,m})\dlog f
    =\frac{1}{f}(u^{e-1-M}\omega_1\wedge \d u \wedge \d c
    +u^{e-M}\omega_2\wedge \d c
    +ecu^{e-1-M}\omega_2\wedge \d u).
    \]
    Since $e-1-M\geq M-1+\sum_{j=1}^m |m_j'|\geq 0$, $\pi^\ast(\eta_{n,m})\dlog f$ is regular at $x$.
    Hence, from \eqref{eq:res} and the residue theorem, we obtain that
    \begin{align*}
        (i_0^\ast-i_\infty^\ast)\alpha^\ast \eta_{n,m}
        &=\sum_{x\in C\smallsetminus |D_C|}\mathrm{Res}_{C,x}(\pi^\ast(\eta_{n,m})\dlog f )\\
        &=\sum_{x\in C}\mathrm{Res}_{C,x}(\pi^\ast(\eta_{n,m})\dlog f )=0.
    \end{align*}
    This completes the proof.
\end{proof}

By \cite[Lemma~5.1.1]{KMSY3}, the graph of the multiplication map
\[
 \mu^\circ\colon \bA^1\times\bA^1\to\bA^1,
 \qquad
 (x,y)\mapsto xy,
\]
gives a morphism $\Gapair\otimes \Gapair\to \Gapair$ which induces a morphism 
\begin{align}\label{mult-Ga}
    \mu\colon \GaMod\otimes\GaMod\to\GaMod.
\end{align}
For any $n>0$, we have a morphism
\[
\mu\otimes \Id^{\otimes (n-1)}
\colon
(\GaMod)^{\otimes (n+1)}=(\GaMod\otimes \GaMod)\otimes (\GaMod)^{\otimes (n-1)}\to \GaMod\otimes (\GaMod)^{\otimes (n-1)}=(\GaMod)^{\otimes n}.
\]
For any $0\le q<n$, we have the following composite morphism 
\[
\mu_q\colon=(\mu\otimes \Id^{\otimes q})\circ \cdots \circ (\mu\otimes \Id^{\otimes (n-1)})\colon (\GaMod)^{\otimes (n+1)}\to  (\GaMod)^{\otimes (q+1)}.
\]
Put $[n]:=\set{1,\ldots,n}$ for $n\geq1$ and $[0]:=\varnothing$.  
For any $S=\{i_1<\cdots <i_q\}\subset[n]$, we define a morphism 
\[
\mu_S:(\GaMod)^{\otimes (n+1)}\to (\GaMod)^{\otimes (q+1)}
\]
to be the composite
\begin{align}
    (\GaMod)^{\otimes (n+1)}\xrightarrow{f}(\GaMod)^{\otimes (n+1)} \xrightarrow{\mu_q} (\GaMod)^{\otimes (q+1)},
\end{align}
where $f$ is an isomorphism induced by the following permutation map
\[
(\bA^1)^{n+1}\to (\bA^1)^{n+1},\ (t_0, t_1, \dots,t_n)\mapsto (t_0, \boldsymbol{x}, t_{i_1},\dots, t_{i_q}) \quad \boldsymbol{x}=(t_j)_{j\in[n]\smallsetminus S}.
\]
The morphism $\mu_S$ is actually induced by a map
\[
    \mu_S^\circ \colon (\bA^1)^{n+1}\to (\bA^1)^{q+1},\ (t_0, t_1, \dots,t_n)\mapsto (t_0\prod_{j\in[n]\smallsetminus S}t_j, t_{i_1},\dots, t_{i_q}).
    \]
We now define a morphism
\[
\theta_{S,m} \colon \Acal_{n+1,m}\to \Omega^{|S|+m}
\]
to be the composite
\[
\Acal_{n+1,m}\xrightarrow{h_{0,\Nis}(\mu_S\otimes \Id_{\GmMod}^{\otimes m})}\Acal_{q+1,m} \xrightarrow{\theta_{q,m}} \Omega^{|S|+m}.
\]

For a finitely generated field $F/k$, a finite extension $E/F$, 
$a_0,\ldots,a_n\in E$ and $b_1,\dots, b_m\in E^\times$, its value on the bracket of \autoref{thm:pres} is
\begin{equation}\label{eq:thetaF}
 \theta_{S,m,F}([a_0,a_1,\ldots,a_n;b_1,\dots,b_m]_{E/F})
 =
 \Tr_{E/F}\left(
 a_0\prod_{j\in[n]\smallsetminus S}a_j
 \bigwedge_{i\in S}\d a_i
 \bigwedge_{j=1}^m \dlog b_j
 \right).
\end{equation}

Put
\begin{equation}\label{eq:Dsum}
 \Dcal_{n+1,m}:=\bigoplus_{S\subset[n]}\Omega^{|S|+m}
\end{equation}
and define
\begin{equation}\label{eq:Theta}
\Theta_{n+1,m}:=\bigoplus_{S\subset[n]}\theta_{S,m}:
 \Acal_{n+1,m}\to\Dcal_{n+1,m}.
\end{equation}

\subsection{Top-degree inverse maps}\label{sec:inv}

Let $0\leq q\leq n$, $m\ge 0$ and let
$S=\set{i_1<\cdots<i_q}\subset[n]$. 
From \autoref{lem:fixed-insertion}, we have a homomorphism 
\begin{equation}\label{eq:insert}
 \iota^{(n)}_{S,m,F}:\Acal_{q+1,m}(F)\to\Acal_{n+1,m}(F)
\end{equation}
such that
\begin{equation}\label{eq:insfld}
 \iota^{(n)}_{S,m,F}
 ([a_0,a_{i_1},\ldots,a_{i_q};b_1,\dots,b_m]_{E/F})
 :=[a_0,a_1(S),\ldots,a_n(S);b_1,\dots,b_m]_{E/F},
\end{equation}
where 
\[
a_j(S) = \begin{cases}
	a_j & j \in S,\\
	1, & j\notin S,
\end{cases}\qquad (j\in [n]).
\]
For $q=0$, this sends
$[a_0;b_1,\ldots,b_m]_{E/F}$ to
$[a_0,1,\ldots,1;b_1,\ldots,b_m]_{E/F}$.

\begin{lem}\label{lem:insert}
 The map \eqref{eq:insert} is well-defined.  
 It is natural for field
 extensions and commutes with finite transfers.
\end{lem}

\begin{proof}
We prove that the map \eqref{eq:insert} comes from a morphism of reciprocity sheaves
\begin{align}\label{eq:iota}
        \iota_{S,m}^{(n)}: \Acal_{q+1,m}\to \Acal_{n+1,m}.
\end{align}
	Consider the map 
    $$\tilde{\iota}_S^{(n)}: (\bA^1)^{q+1}\to (\bA^1)^{n+1}, (t_0, t_1, \dots, t_q) \mapsto (y_0, y_1, \dots, y_n)$$
    where $y_0 = t_0$, $y_{i_r}=t_r$ for $1\leq r\leq q$ and $y_j=1$ for $j\notin S$.
    The map $\tilde{\iota}_S^{(n)}$ gives a morphism $\Gapair^{\otimes (q+1)}\to \Gapair^{\otimes (n+1)}$.
     By \eqref{eq:splitGaGmMod}, we have the following composite map
    \begin{align}        
    (\GaMod)^{\otimes (q+1)}\otimes(\GmMod)^{\otimes m}
    \to
    \Ztr((\Gapair)^{\otimes (q+1)}&\otimes(\Gmpair)^{\otimes m})\\
    &\xrightarrow{\Ztr(\tilde{\iota}_S^{(n)}\otimes \Id_{\Gmpair}^{\otimes m})}
    \Ztr((\Gapair)^{\otimes (n+1)}\otimes(\Gmpair)^{\otimes m})\\
    &\qquad \qquad\to
    (\GaMod)^{\otimes (n+1)}\otimes(\GmMod)^{\otimes m}.
    \end{align}
    Here the first map is an inclusion of a direct summand, the last map is a projection of a direct summand.
    Now applying $h_{0,\Nis}$ gives the desired map \eqref{eq:iota}.
\end{proof}

\color{black}

For $q=0$ and $m=0$ put
\begin{equation}\label{eq:tau0}
 \sigma_{0,0,F}:F=\Omega_F^0\to\Acal_{1,0}(F),
 \qquad a\mapsto[a]_F.
\end{equation}
For $q=0$ and $m\geq 1$,
\begin{equation}
 \sigma_{0,m,F}:\Omega_F^m\to\Acal_{1,m}(F),
 \qquad a\dlog b_1\wedge\cdots \wedge \dlog b_m\mapsto [a;b_1,\dots,b_m]_F.
\end{equation}
This is the map used in the proof of \autoref{thm:one-additive}.

For $q\geq1$ and $T\subset[q]$, put
\begin{equation}\label{eq:cdef}
 a_i(T)=
 \begin{cases}
 a_i,&i\in T,\\
 1,&i\notin T.
 \end{cases}
\end{equation}
On elementary forms define the top-degree candidate by
\begin{align}\label{eq:tau}
 \sigma_{q,m,F}
 &\left(a\,\bigwedge_{i=1}^q\d a_i\bigwedge_{j=1}^m\dlog b_j\right)\\
 &:=
 \sum_{T\subset[q]}(-1)^{q-|T|}
 \left[
 a\prod_{i\in[q]\smallsetminus T}a_i,
 a_1(T),\ldots,a_q(T);b_1,\dots,b_m
 \right]_F.
\end{align}
When all $a_i$ are nonzero and $A=a a_1\cdots a_q$, write
\begin{equation}\label{eq:deltaq}
 \sigma_{q,m,F}
 \left(A\,\bigwedge_{i=1}^q\dlog a_i\bigwedge_{j=1}^m\dlog b_j\right)
 =\Delta_{q,m,F}(A;a_1,\ldots,a_q;b_1,\dots,b_m),
\end{equation}
where
\begin{align}\label{eq:delq}
 \Delta_{q,m,F}&(A;a_1,\ldots,a_q;b_1,\dots,b_m)\\
 &:=
 \sum_{T\subset[q]}(-1)^{q-|T|}
 \left[
 \frac{A}{\prod_{i\in T}a_i},
 a_1(T),\ldots,a_q(T);b_1,\dots,b_m
 \right]_F.
\end{align}
For $q=1$ and $m=0$, \eqref{eq:tau} is exactly the map $\sigma_{1,0,F}$ of
\eqref{eq:s1}.

\begin{prop}\label{prop:top}
Let $F/k$ be a field. For every $q\geq0$ and $m\geq 0$, the map 
\begin{equation}\label{eq:taumap}
 \sigma_{q,m,F}:\Omega_F^{q+m}\to\Acal_{q+1,m}(F)
\end{equation}
is well-defined and 
an $F$-linear homomorphism. 
It is natural for field extensions and commutes with finite transfers.
\end{prop}
\begin{proof}
For $q=0$ and $m=0$ this is the
identification $\Acal_{1,0}\simeq\Ga=\Omega^0$.
The case $q=0$ and $m\geq 1$ follows from
\autoref{thm:one-additive}, while the case $q=1$ and $m=0$ follows from \autoref{thm:two}.
We henceforth assume $q\ge 1$.

    We use a presentation of $\Omega^{q+m}_F$ in \autoref{app:diff}. 
Put $T_{q+m}(F)=F\otimes_\Z (F^\times)^{\otimes (q+m)}$. 
    \begin{equation}\label{eq:BKmap}
 T_{q+m}(F) \to\Omega_F^{q+m},
 \qquad
 A\otimes c_1\otimes\cdots\otimes c_{q+m}
 \mapsto
 A\,\bigwedge_{i=1}^{q+m}\dlog c_i.
\end{equation}
By \autoref{cor:Aflin}, 
this map is surjective and its
kernel, as an $F$-vector subspace, is generated by the following two
families:
\begin{align}
&A\otimes c_1\otimes\cdots\otimes c_{q+m}
 \quad\text{with }c_i=c_j\text{ for some }i<j,
 \label{eq:A1}\\
&x\otimes x\otimes c_2\otimes\cdots\otimes c_{q+m}
 +(1-x)\otimes(1-x)\otimes c_2\otimes\cdots\otimes c_{q+m},
 \label{eq:A2}
\end{align}
where $x,1-x,c_2,\ldots,c_{q+m}\in F^\times$.
    Consider the map
    \[
    \widetilde \sigma_{q,m,F}: T_{q+m}(F) \to \Acal_{q+1,m}(F)
    \]
    defined by 
    \[
	 A\otimes c_1\otimes\cdots\otimes c_{q+m}
     \mapsto \Delta_{q,m,F}(A;c_1,\ldots,c_q;c_{q+1},\dots,c_{q+m}).
    \]
    
    We first claim that the map $\widetilde \sigma_{q,m,F}$ is a group homomorphism. 
    For $A, A'\in F$ it is clear that the following holds;
    \[
    \widetilde \sigma_{q,m,F}((A+A')\otimes c_1\otimes \cdots \otimes c_{q+m})
    =\widetilde \sigma_{q,m,F}(A\otimes c_1\otimes \cdots \otimes c_{q+m})+\widetilde \sigma_{q,m,F}(A'\otimes c_1\otimes \cdots \otimes c_{q+m}).
    \]
    Similarly, for $C,C'\in F^\times$ it is clear that the following holds;
    \begin{align}
        \widetilde \sigma_{q,m,F}(A\otimes \vc\otimes \cdots \otimes CC'\otimes \cdots \otimes c_{q+m})
        &=\widetilde \sigma_{q,m,F}(A\otimes \vc\otimes \cdots \otimes C\otimes \cdots \otimes c_{q+m})\\
        &\quad +\widetilde \sigma_{q,m,F}(A\otimes \vc\otimes \cdots \otimes C'\otimes \cdots\otimes  c_{q+m}),
    \end{align}
    where $\boldsymbol{c}=c_1\otimes \cdots \otimes c_q$.
    
    Let $x,y\in F^\times$, $1\le i\leq q$ and $c_i=xy$.
    For $i\in T\subset [q]$, put $T'=T\smallsetminus\{i\}$ and  $\displaystyle B_T:=\prod_{r\in T'}c_r$.
    Put $\vc=c_{q+1}\otimes \cdots \otimes c_{q+m}$.
    Then we have
    \begin{align*}
        &\widetilde \sigma_{q,m,F}  (A\otimes  c_1\otimes \cdots \otimes xy\otimes\cdots \otimes \vc)\\
        =&\sum_{i\in T\subset[q]}(-1)^{q-|T|}\Bigl\{
        \left[
        \frac{A}{xyB_T},c_1(T),\ldots, xy,\ldots;\vc
        \right]_F
        -
        \left[
        \frac{A}{B_T},c_1(T'),\ldots,1,\ldots;\vc
        \right]_F \Bigr\}\\
        \stackrel{\eqref{eq:leib}}{=}&\sum_{i\in T\subset[q]}(-1)^{q-|T|}\Bigl\{\left[
        \frac{A}{yB_T},c_1(T),\ldots,y,\ldots;\vc
        \right]_F
        +\left[
        \frac{A}{xB_T},c_1(T),\ldots , x, \ldots;\vc
        \right]_F\\
         &\hspace{80mm}-2\left[
        \frac{A}{B_T},c_1(T),\ldots , 1, \ldots;\vc
        \right]_F \Bigr\}\\
        =&\widetilde \sigma_{q,m,F} (A\otimes  c_1\otimes \cdots \otimes x\otimes \cdots \otimes \vc)
        +\widetilde \sigma_{q,m,F}  (A\otimes  c_1\otimes \cdots \otimes  y\otimes  \cdots \otimes \vc). 
    \end{align*}
    This proves the claim.

    By \autoref{lem:Fvec}, $\widetilde \sigma_{q,m,F}$ is $F$-linear.
    From \autoref{cor:Aflin}, our next task is to show that the map $\widetilde \sigma_{q,m,F}$ sends elements of the form \eqref{eq:A1} and \eqref{eq:A2} to zero.

\smallskip
\noindent
\textbf{(The case $(q,m)=(2,0)$)} 
In the case where $(q,m) = (2,0)$, we have 
\begin{align}\label{eq:s20-tilde}
 \widetilde \sigma_{2,0,F}(A\otimes u\otimes v)
 :={}&[A/(uv),u,v]_F-[A/u,u,1]_F
 -[A/v,1,v]_F+[A,1,1]_F.
\end{align}
For
$x,1-x,c\in F^\times$, we have
\begin{align*}
&\widetilde \sigma_{2,0,F}(x\otimes x\otimes c)
 +\widetilde \sigma_{2,0,F}((1-x)\otimes(1-x)\otimes c)\\
&\quad = [1/c,x,c]_F-[1,x,1]_F-[x/c,1,c]_F+[x,1,1]_F\\
&\quad\quad +[1/c,1-x,c]_F-[1,1-x,1]_F
 -[(1-x)/c,1,c]_F+[1-x,1,1]_F\\
&\quad =[1/c,1,c]_F-[1,1,1]_F-[1/c,1,c]_F+[1,1,1]_F=0.
\end{align*}
Next let $x\in F^\times$ and put $P=A/x^2$.  Then
\[
\widetilde \sigma_{2,0,F}(A\otimes x\otimes x)
 =[P,x,x]_F-[Px,x,1]_F-[Px,1,x]_F+[Px^2,1,1]_F.
\]
By \eqref{eq:sym},
\[
 2[P,x,x]_F=[P,1,x^2]_F+[P,x^2,1]_F.
\]
Therefore
\begin{align*}
2\widetilde \sigma_{2,0,F}(A\otimes x\otimes x)
={}&\bigl([P,1,x^2]_F-2[Px,1,x]_F+[Px^2,1,1]_F\bigr)\\
&+\bigl([P,x^2,1]_F-2[Px,x,1]_F+[Px^2,1,1]_F\bigr)\\
&=0 \qquad \mbox{by \eqref{eq:leib}}.
\end{align*}
Since
$\operatorname{char}(F)\neq2$ and $\Acal_{3,0}(F)$ is an $F$-vector space by
\autoref{lem:Fvec}, we obtain $\widetilde \sigma_{2,0,F}(A\otimes x\otimes x)=0$.  Thus this gives a
well-defined map
\[
 \sigma_{2,0,F}:\Omega_F^2\to\Acal_{3,0}(F).
\]

\smallskip
\noindent
\textbf{(The case $(q,m)=(1,1)$)\ }
In the case where $(q,m) = (1,1)$, we have 
\begin{equation}\label{eq:s11-tilde}
 \widetilde \sigma_{1,1,F}(A\otimes u\otimes v)
 =[A/u,u;v]_F-[A,1;v]_F.
\end{equation}
For the relation \eqref{eq:A2}, the additivity gives
\begin{align*}
&\widetilde \sigma_{1,1,F}(x\otimes x\otimes c)
 +\widetilde \sigma_{1,1,F}((1-x)\otimes(1-x)\otimes c)\\
&=[1,x;c]_F-[x,1;c]_F+[1,1-x;c]_F-[1-x,1;c]_F=0.
\end{align*}
For the relation \eqref{eq:A1}, 
we claim that
\begin{equation}\label{eq:s11-repeated}
 [A/x,x;x]_F-[A,1;x]_F=0
 \qquad(x\in F^\times).
\end{equation}
For $x=1$ this is clear.  
If $x = -1$, 
\[
2([A/(-1),-1;-1]_F - [A,1;-1]_F) = 
[-A,-1;1]_F - [A,1;1]_F = 0.
\] 
Since $2$ is invertible in $F$,
\eqref{eq:s11-repeated} follows.
Assume $x\neq\pm1$ and put
\[
 y:=\frac{x}{1+x}.
\]
Then $y,1-y\in F^\times$ and $y/(1-y)=x$.  Set
\begin{align*}
 M&:= \widetilde \sigma_{1,1,F}(y\otimes y\otimes y)
 +\widetilde \sigma_{1,1,F}((1-y)\otimes(1-y)\otimes(1-y))\\
&= [1,y;y]_F-[y,1;y]_F
 +[1,1-y;1-y]_F-[1-y,1;1-y]_F\\ 
&= (-1)([y,1;y]_F+[1-y,1;1-y]) + ([1,y;y]_F + [1,1-y;1-y]_F) \\
&= 0 \quad (\mbox{by \eqref{eq:diff}}).
\end{align*}
Moreover, additivity gives 
\[
 [1,y;1-y]_F+[1,1-y;1-y]_F=[1,1;1-y]_F
\]
and
\[
 [y,1;1-y]_F+[1-y,1;1-y]_F=[1,1;1-y]_F.
\]
Together with multiplicativity in the last component, these identities give
\begin{equation}\label{eq:s11-M}
 M=[1,y;x]_F-[y,1;x]_F.
\end{equation}
Multiplying \eqref{eq:s11-M} by $(1+x)^2$ and using \eqref{eq:leib},
we obtain
\begin{align*}
0=(1+x)^2M
&=[(1+x)^2,y;x]_F-[x(1+x),1;x]_F\\
&=[1+x,x;x]_F-[x,1+x;x]_F\\
&=[1,x;x]_F-[x,1;x]_F.
\end{align*}
Multiplication by $A/x$ now gives \eqref{eq:s11-repeated}.  Hence
we obtain a well-defined map
\[
 \sigma_{1,1,F}:\Omega_F^2\to\Acal_{2,1}(F).
\]

\smallskip
\noindent
\textbf{(Arbitrary $(q,m)$)}
    We first consider an element of the form \eqref{eq:A2}: $x\otimes x\otimes \vc+ (1-x)\otimes (1-x)\otimes \vc$ with $x, 1-x\in F^\times$ where $\vc=c_2\otimes \cdots \otimes c_{q+m}$.
For $T\subset[q]$ and $2\leq i\leq q$, put
\[
c_i(T)=
\begin{cases}
c_i,&i\in T,\\
1,&i\notin T.
\end{cases}
\]
For simplicity, write
\[
[a,b,\mathbf c(T)]_F
:=
[a,b,c_2(T),\ldots,c_q(T);
 c_{q+1},\ldots,c_{q+m}]_F.
\] 
    For $1\in T\subset [q]$, put
    $
    T':=T\smallsetminus\{1\},\quad 
    B_T:=\prod_{r\in T'}c_r.
    $
    Then we have 
    \begin{align}\label{eq:zero-A2}
        \widetilde \sigma_{q,m,F} & (x\otimes  x\otimes\vc +(1-x)\otimes(1-x)\otimes \vc)\\
        &=\sum_{\{1\}\subset T\subset [q]}\Bigl\{(-1)^{q-|T|}
        \left[
        \frac{1}{B_T},x,\vc(T)
        \right]_F
        +(-1)^{q-|T|}
        \left[
        \frac{1}{B_T},1-x,\vc(T)
        \right]_F \\
        &\qquad\qquad+(-1)^{q-|T'|}
        \left[
        \frac{x}{B_T},1,\vc(T)
        \right]_F 
        +(-1)^{q-|T'|}
        \left[
        \frac{1-x}{B_T},1,\vc(T)
        \right]_F \Bigr\}\\
        &=\sum_{\{1\}\subset T\subset[q]}(-1)^{q-|T|}\Bigl\{
        \left[
        \frac{1}{B_T},1,\vc(T)
        \right]_F
        -
        \left[
        \frac{1}{B_T},1,\vc(T)
        \right]_F \Bigr\}\\
        &=0.
    \end{align}
    
    Next we consider elements of the form \eqref{eq:A1}: $A\otimes c_1\otimes\cdots\otimes c_{q+m}\in T_{q+m}(F)$ with $c_i=c_j=:x$ for some $i<j$.
    
    In case $q+1\leq i$, by \eqref{eq:diff2} in \autoref{lem:tworel}, we have
    \[
    \widetilde \sigma_{q,m,F}(A\otimes c_1\otimes\cdots\otimes c_{q+m})=0.
    \]

In case $j \leq q$, 
for $R\subset[q]\smallsetminus\set{i,j}$ put
$B_R=\prod_{r\in R}c_r$.  
For $r\in[q]\smallsetminus\set{i,j}$, put
\[
 c_r(R):=
 \begin{cases}
 c_r,&r\in R,\\
 1,&r\notin R.
 \end{cases}
\]
We consider, in the defining sum \eqref{eq:delq}, the four terms
indexed by $R, R\cup \set{i}, R\cup \set{j}$ and $R\cup \set{i,j}$. 
The term indexed by $T=R$ is
\[
 (-1)^{q-|R|}
 \left[
 \frac{A}{B_R},
 c_1(R),\ldots,\overset{i}{1},\ldots,\overset{j}{1},\ldots,c_q(R);
 \mathbf b
 \right]_F,
\]
where the entries $1$ displayed above are in the $i$-th and $j$-th
positions, respectively.
The term indexed by $T=R\cup\set{i}$ is
\[
 (-1)^{q-|R|-1}
 \left[
 \frac{A}{xB_R},
 c_1(R),\ldots,\overset{i}x,\ldots,\overset{j}1,\ldots,c_q(R);
 \mathbf b
 \right]_F.
\]
Similarly, the term indexed by $T=R\cup\set{j}$ is
\[
 (-1)^{q-|R|-1}
 \left[
 \frac{A}{xB_R},
 c_1(R),\ldots,\overset{i}{1},\ldots,\overset{j}{x},\ldots,c_q(R);
 \mathbf b
 \right]_F,
\]
and the term indexed by $T=R\cup\set{i,j}$ is
\[
 (-1)^{q-|R|-2}
 \left[
 \frac{A}{x^2B_R},
 c_1(R),\ldots,\overset{i}x,\ldots,\overset{j}x,\ldots,c_q(R);
 \mathbf b
 \right]_F.
\]
The sum of these four terms is
\begin{align*}
(-1)^{q-|R|}
\Biggl\{&
 \left[
 \frac{A}{x^2B_R},
 c_1(R),\ldots,x,\ldots,x,\ldots,c_q(R);
 \mathbf b
 \right]_F\\
&-
 \left[
 \frac{A}{xB_R},
 c_1(R),\ldots,x,\ldots,1,\ldots,c_q(R);
 \mathbf b
 \right]_F\\
&-
 \left[
 \frac{A}{xB_R},
 c_1(R),\ldots,1,\ldots,x,\ldots,c_q(R);
 \mathbf b
 \right]_F\\
&+
 \left[
 \frac{A}{B_R},
 c_1(R),\ldots,1,\ldots,1,\ldots,c_q(R);
 \mathbf b
 \right]_F
\Biggr\}.
\end{align*}
Let $\iota_R:\Acal_{3,0}(F)\to\Acal_{q+1,m}(F)$
be the insertion homomorphism of
\autoref{lem:fixed-insertion} given by 
\[
 \iota_R([u,v,w]_F)
 =
 [u,c_1(R),\ldots,\overset{i}{v},\ldots,\overset{j}{w},\ldots,c_q(R);\mathbf b]_F.
\]
Therefore the above four-term sum is
\begin{align*}
&(-1)^{q-|R|}
\iota_R\Biggl(
 \left[\frac{A}{x^2B_R},x,x\right]_F
 -
 \left[\frac{A}{xB_R},x,1\right]_F
-
 \left[\frac{A}{xB_R},1,x\right]_F
 +
 \left[\frac{A}{B_R},1,1\right]_F
\Biggr) \\
&= (-1)^{q-|R|}\iota_R \widetilde \sigma_{2,0,F}
 \left(
 \frac{A}{B_R}\otimes x\otimes x
 \right) = 0, 
\end{align*}
by the well-definedness of $\sigma_{2,0,F}$.

    Finally, suppose that $i\leq q<j$. 
    For $R\subset[q]\smallsetminus\set{i}$ put
	$B_R=\prod_{r\in R}c_r$.  
	For $r\in[q]\smallsetminus\set{i}$, put
\[
 c_r(R):=
 \begin{cases}
 c_r,&r\in R,\\
 1,&r\notin R.
 \end{cases}
\]
We group together, in the sum \eqref{eq:delq}, the two terms
indexed by
 $R$ and $R\cup \set{i}$.
The term indexed by $T=R$ is
\[
 (-1)^{q-|R|}
 \left[
 \frac{A}{B_R},
 c_1(R),\ldots,\overset{i}{1},\ldots,c_q(R);
 c_{q+1},\ldots,\overset{j}{x},\ldots,c_{q+m}
 \right]_F,
\]
where the displayed $1$ is in the $i$-th position and the
displayed $x$ is in the $j$-th position.
On the other hand, the term indexed by
$T=R\cup\set{i}$ is
\[
 (-1)^{q-|R|-1}
 \left[
 \frac{A}{xB_R},
 c_1(R),\ldots,\overset{i}{x},\ldots,c_q(R);
 c_{q+1},\ldots,\overset{j}x,\ldots,c_{q+m}
 \right]_F,
\]
because $c_i=x$.
Let $\iota_R:\Acal_{2,1}(F)\to\Acal_{q+1,m}(F)$
be the insertion homomorphism of
\autoref{lem:fixed-insertion} defined by
\[
 \iota_R([u,v;w]_F)
 =[u,c_1(R),\ldots,\overset{i}{v},\ldots,c_q(R); c_{q+1},\ldots,\overset{j}{w},\ldots,c_{q+m}]_F.
\]	
Therefore the sum of the above two terms is equal to
\begin{align*}
&(-1)^{q-|R|-1}
\iota_R
\left(
 \left[\frac{A}{xB_R},x;x\right]_F
 -
 \left[\frac{A}{B_R},1;x\right]_F
\right) \\
&\quad =(-1)^{q-|R|-1} \iota_R\widetilde \sigma_{1,1,F}
 \left(
 \frac{A}{B_R}\otimes x\otimes x
 \right) =0,
\end{align*}
by the well-definedness of $\sigma_{1,1,F}$. 

    By the above computation, $\widetilde \sigma_{q,m,F}$ induces an $F$-linear homomorphism
    \[
    \sigma_{q,m,F}:\Omega^{q+m}_F\to \Acal_{q+1,m}(F).
    \]

    Finally, we will show that for any finite field extension $E/F$, the following diagrams commute;
    \begin{align}\label{eq:diag}     
    \xymatrix{
    \Omega^{q+m}_E \ar[d]_{\Tr_{E/F}}\ar[r]^-{\sigma_{q,m,E}} &\Acal_{q+1,m}(E) \ar[d]^{\Tr_{E/F}} && \Omega^{q+m}_E \ar[r]^-{\sigma_{q,m,E}} &\Acal_{q+1,m}(E) \\
    \Omega^{q+m}_F \ar[r]^-{\sigma_{q,m,F}} &\Acal_{q+1,m}(F),&&\Omega^{q+m}_F \ar[u]^{\res_{E/F}}\ar[r]^-{\sigma_{q,m,F}} &\Acal_{q+1,m}(F)\ar[u]_{\res_{E/F}}.
    }
    \end{align}

    For an element $\omega=a\bigwedge_{i=1}^q\d a_i\bigwedge_{j=1}^m\dlog b_j\in \Omega^{q+m}_F$, we have
    \begin{align*}        
    \sigma_{q,m,E}\circ \res_{E/F}(\omega)
    &=\sum_{T\subset [q]}(-1)^{q-|T|}
    \left[
    a\prod_{i\in [q]\smallsetminus T}a_i,
    a_1(T),\ldots,a_q(T);b_1,\dots, b_m\right]_E\\
    &=\res_{E/F}\circ \sigma_{q,m,F}(\omega),
    \end{align*}
    which shows that the right diagram of \eqref{eq:diag} is commutative.

    Next we prove the commutativity of the left diagram of \eqref{eq:diag}.
    First, consider the case where $E/F$ is separable. Since $\Omega_E^1 \simeq E\otimes_F \Omega_F^1$, 
    $\Omega^{q+m}_E$ is generated by elements of the form $\omega=a\bigwedge_{i=1}^q\d a_i\bigwedge_{j=1}^m\dlog b_j$ with $a\in E$, $a_1,\dots, a_q\in F$ and $b_1,\dots,b_m\in F^\times$.
    Thus we have
    \begin{align*}
        \sigma_{q,m,F}\circ\Tr_{E/F}(\omega)
        &=\sigma_{q,m,F}(\Tr_{E/F}(a)\bigwedge_{i=1}^q\d a_i\bigwedge_{j=1}^m\dlog b_j)\\
        &=\sum_{T\subset[q]}(-1)^{q-|T|}
        \left[
        \Tr_{E/F}(a)A_T
        ,a_1(T),\ldots,a_q(T);\vb\right]_F\\
        &=\sum_{T\subset[q]}(-1)^{q-|T|}
        \Tr_{E/F}\left[
        aA_T
        ,a_1(T),\ldots,a_q(T);\vb\right]_E\\
        &=\Tr_{E/F}\circ \sigma_{q,m,E}(\omega).
    \end{align*}
     Here $A_T=\prod_{i\in [q]\smallsetminus T}a_i$.
    
    In case where $E/F$ is purely inseparable of degree $p$, $E$ can be written as $E=F[x]$ with $y=x^p\in F$.
    Thus, for $q\ge 1$, every element of $\Omega_E^{q+m}$ is written as a sum of elements of the form
    \[
    \omega=ax^i\bigwedge_{i=1}^q\d a_i\bigwedge_{j=1}^m\dlog b_j,\quad \omega'=ax^i\d x\bigwedge_{i=2}^q\d a_i\bigwedge_{j=1}^m\dlog b_j
    \]
    where $a, a_1,\dots, a_q\in F$, $b_1,\dots,b_m\in F^\times$ and $i\in\set{0,\dots, p-1}$.
    The above case works for $\omega$.
    For $\omega'$, if $0\le i\le p-2$ then we have
    \[
    \omega'
    =\frac{a}{i+1}\d x^{i+1}\wedge \d\va\wedge\dlog \vb
    \]
    where $\d\va=\bigwedge_{i=2}^q\d a_i, \dlog\vb=\bigwedge_{j=1}^m\dlog b_j.$
    Thus we have
    \[
    \Tr_{E/F}(\omega')=\frac{a}{i+1}\d \Tr_{E/F}(x^{i+1})\d\va\wedge\dlog \vb.
    \]
    Therefore, we have
    \begin{align*}
        \sigma_{q,m,F}\circ\Tr_{E/F}(\omega')
        &=\sigma_{q,m,F}\left(\frac{a}{i+1}\d\Tr_{E/F}(x^{i+1})\d\va\wedge\dlog \vb\right)\\
        &=\sum_{\{1\}\subset T\subset[q]}\Bigl\{(-1)^{q-|T|}
        \left[
        \frac{a}{i+1}A_T
        ,\Tr_{E/F}(x^{i+1}),\va(T);\vb
        \right]_F\\
        &\qquad \qquad \qquad +(-1)^{q-|T'|}
        \left[
        \frac{a\Tr_{E/F}(x^{i+1})}{(i+1)}A_T
        ,1,\va(T);\vb
        \right]_F 
         \Bigr\}\\
        &=\sum_{\{1\}\subset T\subset[q]}\Tr_{E/F}\Bigl\{(-1)^{q-|T|}
        \left[
        \frac{a}{i+1}A_T
        ,x^{i+1},\va(T);\vb
        \right]_E\\
        &\qquad \qquad \qquad +(-1)^{q-|T'|}
        \left[
        \frac{ax^{i+1}}{(i+1)}A_T
        ,1,\va(T);\vb
        \right]_E 
         \Bigr\}\\
         &=\Tr_{E/F}\circ \sigma_{q,m,E}(\omega').
    \end{align*}
    Here $T'=T\smallsetminus\set{1}$, $\va(T):=(a_2(T),\ldots,a_q(T))$ and $A_T=\prod_{i\in [q]\smallsetminus T}a_i$.
    
    Lastly, we consider the element $ax^{p-1}\d x\wedge \d \va\wedge\dlog\vb$.
    We have 
    \[
    \Tr_{E/F}(ax^{p-1}\d x\wedge \d \va\wedge\dlog\vb)=a\d y\wedge \d \va\wedge\dlog\vb
    \]
    and hence
    \begin{align}
        &\sigma_{q,m,F}\circ \Tr_{E/F}(ax^{p-1}\d x\wedge \d \va\wedge\dlog\vb)
        \\&=\sigma_{q,m,F}(a\d y\wedge \d \va\wedge\dlog\vb)\\
        &=\sum_{\{1\}\subset T\subset[q]}(-1)^{q-|T|}\Bigl\{
        \left[
        aA_T
        ,y,\va(T);\vb
        \right]_F
        -
        \left[
        ayA_T
        ,1,\va(T);\vb
        \right]_F 
         \Bigr\}.
    \end{align}
    Here $\va(T):=(a_2(T),\ldots,a_q(T))$ and $A_T=\prod_{i\in [q]\smallsetminus T}a_i$. 
    On the other hand, we have
    \begin{align}
        &\Tr_{E/F}\circ  \sigma_{q,m,E}(ax^{p-1}\d x \wedge \d\va \wedge \dlog \vb)\\
        =&\Tr_{E/F}\sum_{\{1\}\subset T\subset[q]}(-1)^{q-|T|}\Bigl\{
        \left[
        ax^{p-1}A_T
        ,x,\va(T);\vb
        \right]_E 
        -
        \left[
        ayA_T
        ,1,\va(T);\vb
        \right]_E 
         \Bigr\}\\
        =&\sum_{\set{1}\subset T\subset[q]}(-1)^{q-|T|}\Bigl\{
        \Tr_{E/F}\!\left[
        ax^{p-1}A_T
        ,x,\va(T);\vb
        \right]_E\!\!\! 
        -
        \Tr_{E/F}\!\left[
        ayA_T
        ,1,\va(T);\vb
        \right]_E 
         \Bigr\}\\
        =&\sum_{\{1\}\subset T\subset[q]}(-1)^{q-|T|}
        \left[
        ax^{p-1}A_T
        ,x,\va(T);\vb
        \right]_{E/F}.
    \end{align}
    Therefore, it suffices to show that
    \begin{align}\label{eq:rel-insep}        
    [x^{p-1},x,\va;\vb]_{E/F}-[1,y,\va;\vb]_F+[y,1,\va;\vb]_F=0.
    \end{align}
    Consider a homomorphism \eqref{eq:insert-two-Ar}
    \[
    \iota_{1,2}:\Acal_{2,0}(F)\to \Acal_{q+1,m}(F),\quad [u,v]_{L/F} \mapsto
    [u,v,\va;\vb]_{L/F}.
    \]
    Then the left hand side of \eqref{eq:rel-insep} is the image of $\xi:=[x^{p-1},x]_{E/F}-[1,y]_F+[y,1]_F\in \Acal_{2,0}(F)$ via $\iota_{1,2}$.
    By \autoref{thm:two} (cf.~\cite[Proof of Proposition~5.13]{RSY22}), the element $\xi=0$, and hence we have \eqref{eq:rel-insep}.

    For a general finite extension $E/F$, let $E_s/F$ be the maximal
separable subextension.  Then $E/E_s$ is purely inseparable and can
be written as a tower of purely inseparable extensions of degree $p$.
The result follows from the transitivity of the trace map 
and of the transfer on $\Acal_{q+1,m}$.
    This completes the proof.
\end{proof}

We now obtain every component of the inverse from a top-degree map.  If
$S=\set{i_1<\cdots<i_q}\subset[n]$, define
\begin{equation}\label{eq:sigma}
 \sigma_{S,m,F}:=
 \iota^{(n)}_{S,m,F}\circ \sigma_{q,m,F}:
 \Omega_F^{q+m}\to\Acal_{n+1,m}(F).
\end{equation}
Equivalently, if the differential variables are indexed by the elements of
$S$, then
\begin{equation}\label{eq:sigexp}
 \sigma_{S,m,F}
 \left(a\,\bigwedge_{i\in S}\d a_i\bigwedge_{j=1}^m \dlog b_j\right)
 =\sum_{T\subset S}(-1)^{|S|-|T|}
 \left[
 a\prod_{i\in S\smallsetminus T}a_i,
 \va(T);\vb
 \right]_F,
\end{equation}
where 
\[
 \va(T)=(a_i(T))_{i\in[n]},
 \qquad
 a_i(T)=
 \begin{cases}
 a_i,&i\in T,\\
 1,&i\notin T.
 \end{cases}
\]
In particular
$a_j(T)=1$ for $j\notin S$.  For $S=\varnothing$ this gives
\[
 \sigma_{\varnothing,m,F}
 \left(
 a\bigwedge_{j=1}^m\dlog b_j
 \right)
 =
 [a,1,\ldots,1;\vb]_F.
\]

If $a_i\neq0$ for every $i\in S$ and
$A=a\prod_{i\in S}a_i$,
put
\begin{equation}\label{eq:delta}
 \sigma_{S,m,F}
 \left(
 A\bigwedge_{i\in S}\dlog a_i
 \bigwedge_{j=1}^m\dlog b_j
 \right)
 =
 \Delta_{S,m,F}
 \bigl(A;(a_i)_{i\in S};\vb\bigr),
\end{equation}
where
\begin{equation}\label{eq:deldef}
 \Delta_{S,m,F}
 \bigl(A;(a_i)_{i\in S};\vb\bigr)
 :=
 \sum_{T\subset S}(-1)^{|S|-|T|}
 \left[
 \frac{A}{\prod_{i\in T}a_i},
 \va(T);\vb
 \right]_F.
\end{equation}

By \autoref{prop:top} and \autoref{lem:insert}, we obtain the following corollary.
\begin{cor}\label{cor:sigtr}
For every $S\subset[n]$ and $m\geq 0$, the map 
$\sigma_{S,m,F}\colon\Omega_F^{|S|+m}\to\Acal_{n+1,m}(F)$ is well defined and
natural for field extensions. 
Moreover, for every finite extension $E/F$ the diagram
\[
\xymatrix{
 \Omega_E^{|S|+m}\ar[r]^-{\sigma_{S,m,E}}\ar[d]_{\Tr_{E/F}}
 &\Acal_{n+1,m}(E)\ar[d]^{\Tr_{E/F}}\\
 \Omega_F^{|S|+m}\ar[r]^-{\sigma_{S,m,F}}
 &\Acal_{n+1,m}(F)
}
\]
commutes.
\end{cor}
\black

This description makes the induction on the number of additive tensor factors
explicit. 
At the stage of $n+1$ additive factors, if $S\subsetneq[n]$ then
$|S|<n$ and $\sigma_{S,m,F}$ is obtained by inserting units into the
previously constructed top-degree map $\sigma_{|S|,m,F}$.  The only new inverse
component is the one corresponding to $S=[n]$, namely 
\[
 \sigma_{[n],m,F}=\sigma_{n,m,F}:\Omega_F^{n+m}\to\Acal_{n+1,m}(F).
\]

\begin{lem}\label{lem:orth}
For $S,T\subset[n]$, one has
\begin{equation}\label{eq:orth}
 \theta_{T,m,F}\,\sigma_{S,m,F}=
 \begin{cases}
 \Id_{\Omega_F^{|S|+m}},&T=S,\\
 0,&T\neq S.
 \end{cases}
\end{equation}
\end{lem}

\begin{proof}
Write $S=\set{i_1<\cdots<i_q}$ with $q = |S|$ and 
let $\omega:=a\bigwedge_{i \in S}\d a_i\bigwedge_{j=1}^m\dlog b_j\in \Omega_F^{q+m}$, 
where the wedge over $S$ is taken in increasing order. 
By \eqref{eq:sigexp},
\[
\sigma_{S,m,F}(\omega)
=
\sum_{S' \subset S}(-1)^{|S|-|S'|}
\left[a\prod_{i\in S\smallsetminus S'}a_i,\va(S');\vb\right]_F.
\]
From \eqref{eq:thetaF} we obtain
    \begin{align} \label{eq:orth-computation}
      \theta_{T,m,F}\,\sigma_{S,m,F}(\omega)
      &=\theta_{T,m,F}\left(\sum_{S'\subset S}(-1)^{|S|-|S'|}
      \left[
      a\prod_{i\in S\smallsetminus S'}a_i,
      \va(S');\vb
      \right]_F\right)\\
      &=\sum_{S'\subset S}(-1)^{|S|-|S'|}
      \theta_{T,m,F}\left[
      a\prod_{i\in S\smallsetminus S'}a_i,
      \va(S');\vb
      \right]_F\\
      &= \sum_{S'\subset S}(-1)^{|S|-|S'|}
a\prod_{i\in S\smallsetminus S'}a_i
\prod_{j\in[n]\smallsetminus T}a_j(S')
\bigwedge_{i\in T}\d a_i(S')
\bigwedge_{j=1}^m\dlog b_j.
    \end{align}
If $T\not\subset S'$, choose $j\in T\smallsetminus S'$.  Then
$a_j(S')=1$, and therefore
$\bigwedge_{j\in T}\d a_j(S')
=0$.
Thus the $S'$-summand in \eqref{eq:orth-computation} is zero unless
$T\subset S'$.

Suppose first that $T\not\subset S$.  Since $S'\subset S$ for every
summand in \eqref{eq:orth-computation}, we have $T\not\subset S'$ for
every $S'\subset S$.  Hence
$\theta_{T,m,F}\sigma_{S,m,F}(\omega)=0$.

Now assume that $T\subset S$.  By the preceding observation, only the
terms with
$T\subset S'\subset S$
remain.  For such $S'$, we have
$\prod_{j\in[n]\smallsetminus T}a_j(S')
=
\prod_{j\in S'\smallsetminus T}a_j$
and
$\bigwedge_{i\in T}\d a_i(S')
=
\bigwedge_{i\in T}\d a_i$.
Therefore
\begin{align*}
&a\prod_{i\in S\smallsetminus S'}a_i
\prod_{j\in[n]\smallsetminus T}a_j(S')
\bigwedge_{i\in T}\d a_i(S')
\bigwedge_{j=1}^m\dlog b_j\\
&\qquad=
a\left(
\prod_{i\in S\smallsetminus S'}a_i
\right)
\left(
\prod_{j\in S'\smallsetminus T}a_j
\right)
\bigwedge_{i\in T}\d a_i
\bigwedge_{j=1}^m\dlog b_j\\
&\qquad=
a\prod_{i\in S\smallsetminus T}a_i
\bigwedge_{i\in T}\d a_i
\bigwedge_{j=1}^m\dlog b_j.
\end{align*}
Consequently,
\begin{align*}
\theta_{T,m,F}\,\sigma_{S,m,F}(\omega)
&=
\left(
\sum_{T\subset S'\subset S}
(-1)^{|S|-|S'|}
\right)
a\prod_{i\in S\smallsetminus T}a_i
\bigwedge_{i\in T}\d a_i
\bigwedge_{j=1}^m\dlog b_j\\
&=
(1-1)^{|S|-|T|}
a\prod_{i\in S\smallsetminus T}a_i
\bigwedge_{i\in T}\d a_i
\bigwedge_{j=1}^m\dlog b_j.
\end{align*}
If $T\subsetneq S$, then $|S|-|T|>0$, and hence
\[
\theta_{T,m,F}\,\sigma_{S,m,F}(\omega)=0.
\]
If $T=S$, then
\[
\theta_{S,m,F}\,\sigma_{S,m,F}(\omega)
=
a\bigwedge_{i\in S}\d a_i
\bigwedge_{j=1}^m\dlog b_j
=\omega.
\]
Thus
\[
\theta_{T,m,F}\,\sigma_{S,m,F}
=
\begin{cases}
\Id_{\Omega_F^{|S|+m}},&T=S,\\
0,&T\neq S.
\end{cases}
\]
\end{proof}

\subsection{Proof of the main theorem}\label{sec:proof}

\begin{thm}\label{thm:decomp}
Let $k$ be a perfect field of characteristic $p\neq 2$.
Let $n\geq0$ and $m\ge 0$.
There are canonical isomorphisms of Nisnevich reciprocity sheaves
\begin{equation}\label{eq:mainbd}
 h_{0,\Nis}\bigl((\GaSharp)^{\otimes(n+1)}\otimes (\GmSharp)^{\otimes m}\bigr)
 \simeq
 \bigoplus_{S\subset[n]}\Omega^{|S|+m}
 \simeq
 \bigoplus_{q=0}^{n}(\Omega^{q+m})^{\oplus\binom{n}{q}}.
\end{equation}
For a finitely generated field $F/k$, a finite extension
$E/F$, and $S\subset[n]$, the $S$-component is
\begin{equation}\label{eq:fieldbd}
 [a_0,a_1,\ldots,a_n;b_1,\dots, b_m]_{E/F}
 \mapsto
 \Tr_{E/F}\left(
 a_0\prod_{j\in[n]\smallsetminus S}a_j
 \bigwedge_{i\in S}\d a_i
 \bigwedge_{j=1}^m\dlog b_j
 \right),
\end{equation}
where the wedge is taken in increasing order of the indices in $S$.
\end{thm}


\begin{lem}\label{lem:expand}
Let $F/k$ be a field, let $E/F$ be finite, and let
$a_0,a_1,\ldots,a_n\in E$, $b_1,\dots,b_m\in E^\times$.  Assume $a_1,\ldots,a_n\neq0$, and put
$A=a_0a_1\cdots a_n$.  Then in $\Acal_{n+1,m}(F)$,
\begin{equation}\label{eq:expand}
 [a_0,a_1,\ldots,a_n;b_1,\dots,b_m]_{E/F}
 =
 \sum_{S\subset [n]}
 \Delta_{S,m,E/F}(A;(a_i)_{i\in S};\vb).
\end{equation}
Here $\Delta_{S,m,E/F}(A;(a_i)_{i\in S};\vb)=\Tr_{E/F}(\Delta_{S,m,E}(A;(a_i)_{i\in S};\vb))$.
\end{lem}

\begin{proof}
Expand the right-hand side using \eqref{eq:deldef}.  Fix
$T\subset [n]$.  The bracket whose non-unit entries among the last $n$
positions occur exactly in $T$ has coefficient
\begin{equation}\label{eq:altsum2}
 \sum_{T\subset S\subset [n]}(-1)^{|S|-|T|}.
\end{equation}
This coefficient is zero unless $T=[n]$, when it is one.  The unique
surviving term is $[a_0,a_1,\ldots,a_n;b_1,\dots,b_m]_{E/F}$.
\end{proof}

\bigskip

\begin{proof}[Proof of \autoref{thm:decomp}]
The case $n=0$ and $m=0$ is $\Acal_{1,0}\simeq\Ga$,
which follows from the definition of $\GaSharp$; see
\cite[Section~5.1]{RSY22}.  
The case $n=0$ and $m\geq 1$ is \autoref{thm:one-additive}.
We therefore fix an arbitrary integer $n\geq1$.  

Let $F/k$ be a finitely generated field.  For every
$S=\set{i_1<\cdots<i_q}\subset[n]$, \autoref{cor:sigtr} gives a
well-defined $F$-linear map
\[
 \sigma_{S,m,F}:\Omega_F^{q+m}\to\Acal_{n+1,m}(F)
\]
which is natural for field extensions and commutes with finite transfers.
Define
\[
 \Sigma_{n+1,m,F}
 :=
 \bigoplus_{S\subset[n]}\sigma_{S,m,F}:
 \Dcal_{n+1,m}(F)\to\Acal_{n+1,m}(F).
\]
By \autoref{lem:orth}, we have 
\begin{equation}\label{eq:leftn-proof}
 \Theta_{n+1,m,F}\,\Sigma_{n+1,m,F}
 =
 \Id_{\Dcal_{n+1,m}(F)}.
\end{equation}
By \autoref{thm:pres}, the canonical
map
\[
 \bigoplus_{E/F\ {\rm finite}}E^{\otimes_{\Z}(n+1)}\otimes_\Z (E^\times)^{\otimes_{\Z}m}
 \to 
 \Acal_{n+1,m}(F)
\]
is surjective.  Consequently it is enough to prove that for any finite extension $E/F$,
\[
 \Sigma_{n+1,m,F}\,\Theta_{n+1,m,F}([a_0,a_1,\ldots,a_n;b_1,\dots,b_m]_{E/F})=[a_0,a_1,\ldots,a_n;b_1,\dots,b_m]_{E/F}.
\]
If $a_i=0$ for some $0\leq i\leq n$, then the additivity of the bracket
in the $i$-th entry gives
\[
 [a_0,\ldots,0,\ldots,a_n;\vb]_{E/F}
 =
 [a_0,\ldots,0+0,\ldots,a_n;\vb]_{E/F}
 =
 2[a_0,\ldots,0,\ldots,a_n;\vb]_{E/F},
\]
and hence $[a_0,a_1,\ldots,a_n;\vb]_{E/F} =0$. 
The desired identity is then immediate.  
We may therefore assume that
$a_0,\ldots,a_n$ are nonzero.  Put
\[
 A:=a_0a_1\cdots a_n\in E^\times.
\]
For every $S\subset[n]$, the formula
\eqref{eq:thetaF} gives
\begin{align} \label{eq:thlog-proof}
 \theta_{S,m,F}([a_0,a_1,\ldots,a_n;\vb]_{E/F})
 &=
 \Tr_{E/F}\left(
 a_0\prod_{j\in[n]\smallsetminus S}a_j
 \bigwedge_{i\in S}\d a_i\bigwedge_{j=1}^m\dlog b_j
 \right)\\
 &=
 \Tr_{E/F}\left(
 A\bigwedge_{i\in S}\dlog a_i\bigwedge_{j=1}^m\dlog b_j
 \right).
\end{align}
Since $\sigma_{S,m,F}$ commutes with finite transfers by
\autoref{cor:sigtr}, we obtain
\begin{align}
 \sigma_{S,m,F}\,\theta_{S,m,F}([a_0,a_1,\ldots,a_n;\vb]_{E/F})
 &=
 \sigma_{S,m,F}
 \Tr_{E/F}\left(
 A\bigwedge_{i\in S}\dlog a_i\bigwedge_{j=1}^m\dlog b_j
 \right)\notag\\
 &=
 \Tr_{E/F}
 \sigma_{S,m,E}\left(
 A\bigwedge_{i\in S}\dlog a_i\bigwedge_{j=1}^m\dlog b_j
 \right).
 \label{eq:sigtheta-proof}
\end{align}
By \eqref{eq:delta} and \eqref{eq:deldef}, the last term is
\begin{align}
 \Tr_{E/F}
 \sum_{T\subset S}(-1)^{|S|-|T|}
 \left[
 \frac{A}{\prod_{i\in T}a_i},
 a_1(T),\ldots,a_n(T);\vb
 \right]_E
 \notag\\
 =
 \sum_{T\subset S}(-1)^{|S|-|T|}
 \left[
 \frac{A}{\prod_{i\in T}a_i},
 a_1(T),\ldots,a_n(T);\vb
 \right]_{E/F},
 \label{eq:delta-transfer-proof}
\end{align}
where
\[
 a_i(T)=
 \begin{cases}
 a_i,&i\in T,\\
 1,&i\notin T.
 \end{cases}
\]
Thus
\begin{equation}\label{eq:sigtheta-delta-proof}
 \sigma_{S,m,F}\,\theta_{S,m,F}([a_0,a_1,\ldots,a_n;\vb]_{E/F})
 =
 \Delta_{S,m,E/F}(A;(a_i)_{i\in S};\vb).
\end{equation}

Summing \eqref{eq:sigtheta-delta-proof} over all
$S\subset[n]$, we get
\begin{align}
 \Sigma_{n+1,m,F}\,\Theta_{n+1,m,F}([a_0,a_1,\ldots,a_n;\vb]_{E/F})
 &=
 \sum_{S\subset[n]}
 \Delta_{S,m,E/F}(A;(a_i)_{i\in S};\vb)\notag\\
 &=
 [a_0,a_1,\ldots,a_n;\vb]_{E/F},
 \label{eq:rightn-proof}
\end{align}
where the second equality is \autoref{lem:expand}.  Hence
\[
 \Sigma_{n+1,m,F}\,\Theta_{n+1,m,F}
 =
 \Id_{\Acal_{n+1,m}(F)}.
\]
Together with \eqref{eq:leftn-proof}, this proves that
\[
 \Theta_{n+1,m,F}:
 \Acal_{n+1,m}(F)
 \xrightarrow{\simeq}
 \Dcal_{n+1,m}(F)
\]
is an isomorphism for every finitely generated field $F/k$.

Since
\[
 \Theta_{n+1,m}:\Acal_{n+1,m}\to\Dcal_{n+1,m}
\]
is a morphism in $\RSC_{\Nis}$, \cite[Lemma~4.16]{RSY22} implies that
$\Theta_{n+1,m}$ is an isomorphism.  Therefore
\[
 \Acal_{n+1,m}
 \simeq
 \bigoplus_{S\subset[n]}\Omega^{|S|+m}.
\]
There are $\binom nq$ subsets $S\subset[n]$ of cardinality $q$.
Grouping the summands according to $|S|$ gives
\[
 \Acal_{n+1,m}
 \simeq
 \bigoplus_{q=0}^{n}
 (\Omega^{q+m})^{\oplus\binom nq}.
\]
\end{proof}

\section{Hochschild homology}\label{sec:HH}

We apply \autoref{thm:decomp} to the Hochschild chain complex determined by
multiplication on $\GaSharp$.  
We use the convention of \cite[Section~1.1]{Lod98}.
Recall that $\CI$ is a Serre subcategory of
$\MPST$, hence an abelian category; see \cite[Section~3.6]{RSY22}.  Thus
homology objects of complexes in $\CI$ are defined. 

We first recall the multiplication on $\GaSharp$.
The multiplication map
\[
 \mu^\circ:\bA^1\times\bA^1\to\bA^1,
 \qquad
 (x,y)\mapsto xy,
\]
induces a
morphism $\GaMod\otimes\GaMod\to\GaMod$ (cf.~\eqref{mult-Ga}).
Applying $h_0^{\bcube}$ gives
\begin{equation}\label{eq:HH-mult}
 \mu:\GaSharp\overset{\CI}{\otimes}\GaSharp\to\GaSharp.
\end{equation}
Associativity follows from the associativity of multiplication on
$\bA^1$.

\subsection{Hochschild complexes}
We regard
$\GaSharp\otimes(\GmSharp)^{\otimes m}$
as a $\GaSharp$-bimodule, with both actions induced
by $\mu$ on the first factor.
Let $n\geq0$ and $m\geq0$, and put
\begin{equation}\label{eq:HH-chain}
 C_n^{(m)}(\GaSharp) :=
 (\GaSharp)^{\overset{\CI}{\otimes} (n+1)}\overset{\CI}{\otimes}(\GmSharp)^{\overset{\CI}{\otimes} m}.
\end{equation}
For $n\geq1$ and $0\leq i<n$, define
\begin{equation}\label{eq:HH-face}
 d_i
 :=
 \Id_{\GaSharp}^{\otimes i}
 \otimes\mu
 \otimes
 \Id_{\GaSharp}^{\otimes(n-i-1)}
 \otimes
 \Id_{\GmSharp}^{\otimes m}:
 C_n^{(m)}(\GaSharp)
 \to
 C_{n-1}^{(m)}(\GaSharp).
\end{equation}
For the last face $d_n$,
let
\[
 \tau_n:
 (\GaSharp)^{\overset{\CI}{\otimes}(n+1)}
 \overset{\CI}{\otimes}(\GmSharp)^{\overset{\CI}{\otimes} m}
 \to
 (\GaSharp)^{\overset{\CI}{\otimes}(n+1)}
 \overset{\CI}{\otimes}(\GmSharp)^{\overset{\CI}{\otimes} m}
\]
be the symmetry isomorphism given by
\begin{equation}\label{eq:HH-cyclic}
 \tau_n
 (a_0\otimes a_1\otimes\cdots\otimes a_n\otimes\vb)
 =
 a_n\otimes a_0\otimes a_1\otimes\cdots\otimes a_{n-1}\otimes\vb.
\end{equation}
Define
$d_n
 :=
 \left(
 \mu\otimes
 \Id_{\GaSharp}^{\otimes(n-1)}
 \otimes
 \Id_{\GmSharp}^{\otimes m}
 \right)
 \circ\tau_n$.
For $n\geq1$, the boundary is 
\begin{equation}\label{eq:HH-boundary}
 \partial_n=\sum_{i=0}^{n}(-1)^i d_i:
 C_n^{(m)}(\GaSharp)\to C_{n-1}^{(m)}(\GaSharp).
\end{equation}
We set $\partial_0 = 0$.
Thus, after applying
$a_{\Nis}\omega_{\CI}$, 
for a finitely
generated field $F/k$ and a finite extension $E/F$, one has
\begin{align}
 d_{0,F}
 [a_0,a_1,\ldots,a_n;\vb]_{E/F}
 &=
 [a_0a_1,a_2,\ldots,a_n;\vb]_{E/F},
 \label{eq:face-zero}\\
 d_{i,F}
 [a_0,a_1,\ldots,a_n;\vb]_{E/F}
 &=
 [a_0,\ldots,a_{i-1},a_ia_{i+1},
 a_{i+2},\ldots,a_n;\vb]_{E/F}
 \quad (1\leq i<n),
 \label{eq:face-middle}\\
 d_{n,F}
 [a_0,a_1,\ldots,a_n;\vb]_{E/F}
 &=
 [a_na_0,a_1,\ldots,a_{n-1};\vb]_{E/F}.
 \label{eq:face-last}
\end{align}
Hence
\begin{align}
 \partial_{n,F}
 [a_0,a_1,\ldots,a_n;\vb]_{E/F}
 \notag
 &=
 [a_0a_1,a_2,\ldots,a_n;\vb]_{E/F}
 \notag\\
 &\quad+
 \sum_{i=1}^{n-1}(-1)^i
 [a_0,\ldots,a_{i-1},a_ia_{i+1},
 a_{i+2},\ldots,a_n;\vb]_{E/F}.
  \notag\\
  &\quad+
  (-1)^n
  [a_na_0,a_1,\ldots,a_{n-1};\vb]_{E/F}. \label{eq:HH-boundary-field}
\end{align}

We denote the $n$-th homology object of this complex in $\CI$ by
\begin{equation}\label{eq:HHm-def}
 \mathrm{HH}_n^{(m)}(\GaSharp)
 :=
 H_n\bigl(C_\bullet^{(m)}(\GaSharp)\bigr).
\end{equation}
 Notice the index shift
\begin{equation}\label{eq:HH-shift}
 a_{\Nis}\omega_{\CI}C_n^{(m)}(\GaSharp)=\Acal_{n+1,m} 
\end{equation}
the $n$-th Hochschild chain object has $n+1$ tensor factors of $\GaSharp$.

Put 
\begin{equation}\label{eq:HH-B}
 B:=\Z[\varepsilon]/(\varepsilon^2),
\end{equation} 
and define $\epsilon\colon B\to \Z$ by $\epsilon(1) = 1, \epsilon(\varepsilon) = 0$. 
Consider $B = \Z\oplus \Z\varepsilon$ as a graded $\Z$-algebra by 
grading $\deg(1) = 0$ and $\deg(\varepsilon) = 1$; 
$B_0 = \Z\cdot 1$, $B_1 = \Z\varepsilon$, and $B_q = 0$ for $q\ge 2$. 
If $M$ is an abelian group, 
we regard $M$ as a $B$-bimodule through $\epsilon$, that is,
\[
 b\cdot m=\epsilon(b)m,\qquad m\cdot b=\epsilon(b)m
 \qquad (b\in B,\ m\in M).
\]
In particular, $\varepsilon$ acts trivially on $M$ from both sides.
Thus the Hochschild complex with coefficients in $M$ (\cite[Section~1.1.1--1.1.3]{Lod98}) is 
\begin{equation}\label{eq:HH-Loday}
 C_n(B,M)=M\otimes_{\Z}B^{\otimes n}.
\end{equation}
with boundary
\begin{align}
 b_n(m\otimes e_1\otimes\cdots\otimes e_n)
 ={}&\epsilon(e_1)m\otimes e_2\otimes\cdots\otimes e_n\notag\\
 &+\sum_{i=1}^{n-1}(-1)^i
 m\otimes e_1\otimes\cdots\otimes(e_ie_{i+1})
 \otimes\cdots\otimes e_n\notag\\
 &+(-1)^n\epsilon(e_n)m\otimes
 e_1\otimes\cdots\otimes e_{n-1}.
\end{align}
For $S\subset[n]$, let
$e_S=e_1\otimes\cdots\otimes e_n$, where $e_i=\varepsilon$ if $i\in S$
and $e_i=1$ otherwise.  For $0\leq q\leq n$, put
\begin{equation}\label{eq:HH-Cmq}
 C_{n,q}(B,M)
 :=\bigoplus_{\substack{S\subset[n]\\|S|=q}} M\otimes\Z e_S
 \subset C_n(B,M).
\end{equation}
For $n=0$, we use $C_{0,0}(B,M)=M$ and $C_{0,q}(B,M)=0$ for $q\neq0$.
The complex $C_{\bullet,q}(B,M)$ is a subcomplex of $C_\bullet(B,M)$.

\begin{lem}\label{lem:HHnorm}
Let $M$ be an abelian group, regarded as a $B$-bimodule through $\epsilon$.
For every $q\geq0$,
\begin{equation}\label{eq:HH-norm-hom}
 H_n\bigl(C_{\bullet,q}(B,M)\bigr)
 \simeq
 \begin{cases}
 M,&n=q,\\
 0,&n\neq q.
 \end{cases}
\end{equation}
Under the isomorphism for $n=q$, the class corresponding to $x\in M$ is
represented by $x\otimes\varepsilon^{\otimes q}$.
\end{lem}
\begin{proof}
Recall that
\[
 B=\Z[\varepsilon]/(\varepsilon^2),
 \qquad
 \epsilon(1)=1,\quad \epsilon(\varepsilon)=0.
\]
We use the normalized Hochschild complex of
\cite[Definition~1.1.14]{Lod98}. 
Put 
\[
\overline B:= B/ (\Z\cdot 1) = \Z \varepsilon.
\]
Hence the normalized Hochschild complex with coefficients in $M$ is
\begin{equation}\label{eq:HH-normalized}
 \overline C_n(B,M)
 =M\otimes_{\Z}\overline B^{\otimes n}
 \simeq M\otimes_{\Z}(\Z\varepsilon)^{\otimes n}.
\end{equation}
Its differential is zero. 
Indeed,
$\overline B=\Z\varepsilon$, so $\overline C_n(B,M)$ is generated by
elements of the form
$x\otimes\varepsilon^{\otimes n}$ for $x\in M$.
For such an element, the boundary is 
\begin{align*}
 b_n(x\otimes\varepsilon^{\otimes n})
 ={}&
 x\varepsilon\otimes\varepsilon^{\otimes(n-1)}\\
 &+\sum_{i=1}^{n-1}(-1)^i
 x\otimes\varepsilon^{\otimes(i-1)}
 \otimes\varepsilon^2
 \otimes\varepsilon^{\otimes(n-i-1)}\\
 &+(-1)^n\varepsilon x\otimes\varepsilon^{\otimes(n-1)}.
\end{align*}
Since the left and right actions of $B$ on $M$ are induced by the
$\epsilon$, we have 
\[
 x\varepsilon=\epsilon(\varepsilon)x=0,
 \quad 
 \varepsilon x=\epsilon(\varepsilon)x=0, 
\]
and $\varepsilon^2 = 0$. 

For $q\geq0$, define
\[
 \overline C_{n,q}(B,M)
 :=
 \begin{cases}
  \overline C_n(B,M),&n=q,\\
  0,&n\neq q.
 \end{cases}
\]
Since the differential of the normalized complex is zero, these are
subcomplexes and
\[
 \overline C_\bullet(B,M)
 =
 \bigoplus_{q\geq0}\overline C_{\bullet,q}(B,M).
\]
Moreover,
\[
 \overline C_{n,q}(B,M)
 \simeq
 \begin{cases}
  M,&n=q,\\
  0,&n\neq q,
 \end{cases}
\]
where the isomorphism for $n=q$ sends
$x\otimes\varepsilon^{\otimes q}$ to $x$.

Let
\[
 \pi_n:C_n(B,M)\to\overline C_n(B,M)
\]
be the canonical quotient map induced by $B\to\overline B$; explicitly,
\[
 \pi_n(x\otimes b_1\otimes\cdots\otimes b_n)
 =
 x\otimes\bar b_1\otimes\cdots\otimes\bar b_n, 
\]
where $\bar b_i$ is the image of $b_i$ in $\overline B$. 
These maps form a chain map
\[
 \pi:C_\bullet(B,M)\to\overline C_\bullet(B,M).
\]
By \cite[Proposition~1.1.15]{Lod98}, $\pi$ is a quasi-isomorphism. 
The projection $\pi$ is compatible with the decomposition
\[
 C_\bullet(B,M)
 =
 \bigoplus_{q\geq0}C_{\bullet,q}(B,M).
\]
Therefore $\pi$ restricts to chain maps
\[
 \pi_q:
 C_{\bullet,q}(B,M)
 \to
 \overline C_{\bullet,q}(B,M),
\]
where $\overline C_{\bullet,q}(B,M)$ denotes the corresponding part of
the normalized complex.
Since
\[
 \pi=\bigoplus_{q\geq0}\pi_q
\]
and $\pi$ is a quasi-isomorphism, each $\pi_q$ is a quasi-isomorphism.
Indeed, homology commutes with these direct sums, so for every $n$,
\[
 H_n(\pi)
 =
 \bigoplus_{q\geq0}H_n(\pi_q).
\]
Since $H_n(\pi)$ is an isomorphism, each $H_n(\pi_q)$ is an
isomorphism. 
 Therefore
\[
 H_n\bigl(C_{\bullet,q}(B,M)\bigr)
 \simeq
 H_n\bigl(\overline C_{\bullet,q}(B,M)\bigr)
 \simeq
 \begin{cases}
  M,&n=q,\\
  0,&n\neq q.
 \end{cases}
\]

Finally, when $n=q$, the map $\pi_q$ sends
$x\otimes\varepsilon^{\otimes q}$ 
to the same element in $\overline C_q(B,M)$, which corresponds to
$x\in M$ under the above identification.  Hence the homology class
corresponding to $x$ is represented in $C_{\bullet,q}(B,M)$ by
$x\otimes\varepsilon^{\otimes q}$.
\end{proof}

\subsection{Homology of the decomposed complex}

For a finitely generated field $F/k$, \autoref{thm:decomp} gives
\begin{equation}\label{eq:HH-Theta}
 \Theta_{n+1,m,F}:
 \Acal_{n+1,m}(F)
 \xrightarrow{\simeq}
 \bigoplus_{S\subset[n]}\Omega_F^{|S|+m}.
\end{equation}
We identify the summand indexed by $S$ with
$\Omega_F^{|S|+m}\otimes\Z e_S$.
Thus \eqref{eq:HH-Theta} becomes
\begin{equation}\label{eq:HH-Theta-B}
 \Acal_{n+1,m}(F)
 \xrightarrow{\simeq}
 \bigoplus_{q\geq0}
 C_{n,q}(B,\Omega_F^{q+m}).
\end{equation}

\begin{prop}\label{prop:HHcomplex}
The isomorphisms \eqref{eq:HH-Theta-B} give an isomorphism of complexes
\begin{equation}\label{eq:HH-complex-iso}
 \bigl(\Acal_{\bullet+1,m}(F),\partial_\bullet\bigr)
 \simeq
 \bigoplus_{q\geq0}
 C_{\bullet,q}(B,\Omega_F^{q+m}).
\end{equation}
\end{prop}
\begin{proof}
For $n=0$, both differentials are zero.
For $n\geq1$, we show that the following diagram
is commutative:
\[
\xymatrix{
 \Acal_{n+1,m}(F)
 \ar[r]^-{\partial_{n,F}}
 \ar[d]_{\Theta_{n+1,m,F}}
 &
 \Acal_{n,m}(F)
 \ar[d]^{\Theta_{n,m,F}}
 \\
 \displaystyle\bigoplus_{q\geq0}
 C_{n,q}(B,\Omega_F^{q+m})
 \ar[r]^-{b_n}
 &
 \displaystyle\bigoplus_{q\geq0}
 C_{n-1,q}(B,\Omega_F^{q+m}).
}
\]

Fix $T\subset[n-1]$ and put
$T+1:=\set{j+1\mid j\in T}$.
For $1\leq i\leq n-1$, put
\[
S(T,i) := \set{j\in T\mid j<i}
 \cup
 \set{j+1\mid j\in T,\ j>i}.
\]
Using \eqref{eq:fieldbd}, we obtain 
\begin{equation}\label{eq:HH-face-comp}
 \theta_{T,m,F}\,d_{0,F}
 =
 \theta_{T+1,m,F},
 \qquad
 \theta_{T,m,F}\,d_{n,F}
 =
 \theta_{T,m,F},
\end{equation}
and, for $1\leq i\leq n-1$,
\begin{equation}\label{eq:HH-face-comp-middle}
 \theta_{T,m,F}\,d_{i,F}
 =
 \begin{cases}
\theta_{S(T,i),m,F},&\mbox{if $i\notin T$},\\
  \theta_{S(T,i)\cup\set{i},m,F}
  +
  \theta_{S(T,i)\cup\set{i+1},m,F},
  &\mbox{if $i\in T$}.
 \end{cases}
\end{equation}
Indeed, the first two identities follow directly from
\eqref{eq:face-zero} and \eqref{eq:face-last}.  For the middle faces $1\le i\le n-1$, 
if $i\notin T$, the factor $a_i a_{i+1}$ occurs in the coefficient of
the differential form, and hence
\[
 \theta_{T,m,F}\,d_{i,F}
 =
 \theta_{S(T,i),m,F}.
\]
If $i\in T$, the factor $a_i a_{i+1}$ occurs as a differential, and
the Leibniz rule
$\d(a_i a_{i+1})=a_i\,\d a_{i+1}+a_{i+1}\,\d a_i$ 
gives
\[
 \theta_{T,m,F}\,d_{i,F}
 =
 \theta_{S(T,i)\cup\set{i},m,F}
 +
 \theta_{S(T,i)\cup\set{i+1},m,F}.
\]

We compare these formulas with the differential on
$ \bigoplus_{q\geq0} C_{\bullet,q}(B,\Omega_F^{q+m})$. 
Recall that
\[
 C_{n,q}(B,\Omega_F^{q+m})
 =
 \bigoplus_{\substack{S\subseteq[n]\\|S|=q}}
 \Omega_F^{q+m}\otimes_{\Z}e_S,
\]
where 
\[
 e_S=e_1\otimes\cdots\otimes e_n,
 \qquad
 e_j=
 \begin{cases}
 \varepsilon,&j\in S,\\
 1,&j\notin S.
 \end{cases}
\]
Similarly,
\[
 C_{n-1,q}(B,\Omega_F^{q+m})
 =
 \bigoplus_{\substack{T\subseteq[n-1]\\|T|=q}}
 \Omega_F^{q+m}\otimes_{\Z}e_T.
\]
Hence
\[
 \bigoplus_{q\geq0}
 C_{n-1,q}(B,\Omega_F^{q+m})
 =
 \bigoplus_{T\subseteq[n-1]}
 \Omega_F^{|T|+m}\otimes_{\Z}e_T.
\]
Thus it is enough to compare, for every $T\subseteq[n-1]$, the
components of
\[
 \Theta_{n,m,F}\,\partial_{n,F}(\alpha)
 \quad\mbox{and}\quad
 b_n\Theta_{n+1,m,F}(\alpha)
\]
in the direct summand
$\Omega_F^{|T|+m}\otimes_{\Z}e_T$,
where $\alpha\in\Acal_{n+1,m}(F)$.

First, since
\[
 \partial_{n,F}
 =
 \sum_{i=0}^{n}(-1)^i d_{i,F},
\]
we have
\begin{align*}
\theta_{T,m,F}\,\partial_{n,F}(\alpha)
&=
\theta_{T,m,F}\,d_{0,F}(\alpha)
+
\sum_{i=1}^{n-1}(-1)^i
\theta_{T,m,F}\,d_{i,F}(\alpha)
+(-1)^n\theta_{T,m,F}\,d_{n,F}(\alpha).
\end{align*}
By \eqref{eq:HH-face-comp} and
\eqref{eq:HH-face-comp-middle}, this is
\begin{align}
\theta_{T,m,F}\,\partial_{n,F}(\alpha)
&=
\theta_{T+1,m,F}(\alpha)\notag\\
&\quad+
\sum_{\substack{1\leq i\leq n-1\\i\notin T}}
(-1)^i\theta_{S(T,i),m,F}(\alpha)\notag\\
&\quad+
\sum_{\substack{1\leq i\leq n-1\\i\in T}}
(-1)^i
\Bigl(
\theta_{S(T,i)\cup\set{i},m,F}(\alpha)
+
\theta_{S(T,i)\cup\set{i+1},m,F}(\alpha)
\Bigr)\notag\\
&\quad+
(-1)^n\theta_{T,m,F}(\alpha).
\label{eq:HH-left-T}
\end{align}

We now compute the corresponding component of
$b_n\Theta_{n+1,m,F}(\alpha)$.  By \eqref{eq:HH-Theta-B},
\[
 \Theta_{n+1,m,F}(\alpha)
 =
 \sum_{S\subseteq[n]}
 \theta_{S,m,F}(\alpha)\otimes e_S.
\]
Thus, writing
\[
 e_S=e_1\otimes\cdots\otimes e_n,
\]
the definition of $b_n$ gives
\begin{align}
b_n\Theta_{n+1,m,F}(\alpha)
&=
\sum_{S\subseteq[n]}
\epsilon(e_1)\theta_{S,m,F}(\alpha)
\otimes e_2\otimes\cdots\otimes e_n
\notag\\
&\quad+
\sum_{S\subseteq[n]}
\sum_{i=1}^{n-1}(-1)^i
\theta_{S,m,F}(\alpha)
\otimes
e_1\otimes\cdots\otimes
(e_i e_{i+1})
\otimes\cdots\otimes e_n
\notag\\
&\quad+
(-1)^n
\sum_{S\subseteq[n]}
\epsilon(e_n)\theta_{S,m,F}(\alpha)
\otimes e_1\otimes\cdots\otimes e_{n-1}.
\label{eq:HH-right-expansion}
\end{align}

We determine the component of \eqref{eq:HH-right-expansion} in
$\Omega_F^{|T|+m}\otimes_{\Z}e_T$.
For the first term of $b_n\Theta_{n+1,m,F}(\alpha)$, 
since 
\[
 \epsilon(e_1)=
 \begin{cases}
 1,&1\notin S,\\
 0,&1\in S,
 \end{cases}
\]
a nonzero term has $1\notin S$. 
Moreover 
\[
 e_2\otimes\cdots\otimes e_n=e_T
 \quad\Longleftrightarrow\quad
 S=T+1.
\]
Hence the contribution of the first term is
$ \theta_{T+1,m,F}(\alpha)$.
Similarly, for the last term, we have
\[
 e_1\otimes\cdots\otimes e_{n-1}=e_T
 \quad\Longleftrightarrow\quad
 S=T.
\]
Hence the contribution of the last term is
$(-1)^n\theta_{T,m,F}(\alpha)$. 

It remains to consider the $i$-th term in the middle sum of
\eqref{eq:HH-right-expansion}, where $1\leq i\leq n-1$.
Recall that
\[
 1\cdot1=1,\qquad
 1\cdot\varepsilon=\varepsilon,\qquad
 \varepsilon\cdot1=\varepsilon,\qquad
 \varepsilon^2=0.
\]
Suppose first that $i\notin T$.  Then the $i$-th factor of $e_T$ is
$1$, so we must have
$ e_i=e_{i+1}=1$. 
The remaining factors agree with $e_T$ precisely when $S=S(T,i)$.
Hence the contribution of the $i$-th middle term is
$(-1)^i\theta_{S(T,i),m,F}(\alpha)$. 

Suppose next that $i\in T$.  Then the $i$-th factor of $e_T$ is
$\varepsilon$.  
Thus $ e_i e_{i+1}=\varepsilon$ holds precisely in the two cases
\[
 (e_i,e_{i+1})=(\varepsilon,1)
 \quad\mbox{or}\quad
 (e_i,e_{i+1})=(1,\varepsilon).
\]
These correspond respectively to
$S=S(T,i)\cup\set{i}$ and $S=S(T,i)\cup\set{i+1}$.
Hence the contribution of the $i$-th middle term is
\[
 (-1)^i
 \Bigl(
 \theta_{S(T,i)\cup\set{i},m,F}(\alpha)
 +
 \theta_{S(T,i)\cup\set{i+1},m,F}(\alpha)
 \Bigr).
\]

Combining these, 
we see that the component of
$b_n\Theta_{n+1,m,F}(\alpha)$ in
$\Omega_F^{|T|+m}\otimes_{\Z}e_T$ is
\begin{align*}
&\theta_{T+1,m,F}(\alpha)\\
&\quad+
\sum_{\substack{1\leq i\leq n-1\\i\notin T}}
(-1)^i\theta_{S(T,i),m,F}(\alpha)\\
&\quad+
\sum_{\substack{1\leq i\leq n-1\\i\in T}}
(-1)^i
\Bigl(
\theta_{S(T,i)\cup\set{i},m,F}(\alpha)
+
\theta_{S(T,i)\cup\set{i+1},m,F}(\alpha)
\Bigr)\\
&\quad+
(-1)^n\theta_{T,m,F}(\alpha).
\end{align*}
By \eqref{eq:HH-left-T}, this is exactly
\[
 \theta_{T,m,F}\,\partial_{n,F}(\alpha).
\]
Since this holds for every $T\subseteq[n-1]$, we obtain
\[
 \Theta_{n,m,F}\,\partial_{n,F}
 =
 b_n\Theta_{n+1,m,F}.
\]
Thus the isomorphisms
\eqref{eq:HH-Theta-B} give an isomorphism of complexes
\[
 \bigl(\Acal_{\bullet+1,m}(F),\partial_\bullet\bigr)
 \simeq
 \bigoplus_{q\geq0}
 C_{\bullet,q}(B,\Omega_F^{q+m}).
\]
\end{proof}

\begin{thm}\label{thm:HH}
Let $k$ be a perfect field of characteristic $p\neq 2$. 
For every $n,m\geq0$, there is a canonical isomorphism of Nisnevich
reciprocity sheaves
\begin{equation}\label{eq:HH-main}
 h_{0,\Nis}\bigl(\mathrm{HH}_n^{(m)}(\GaSharp)\bigr)\simeq\Omega^{n+m}.
\end{equation}
On a finitely generated field $F/k$, 
this isomorphism is induced by $\theta_{[n],m,F}$. 
For a finite extension $E/F$, this homomorphism is 
given by 
\begin{equation}\label{eq:HH-HKR-map}
 [a_0,\ldots,a_n;b_1,\ldots,b_m]_{E/F}
 \mapsto
 \Tr_{E/F}(a_0\,\d a_1\wedge\cdots\wedge \d a_n \wedge
 \dlog b_1\wedge\cdots\wedge\dlog b_m).
\end{equation}
\end{thm}

\begin{proof}
By \cite[Section~3.11, (3.11.1)]{RSY22},
$\omega_{\CI}\colon\CI\to\RSC$ is exact.  Since Nisnevich sheafification of abelian
presheaves is exact, the composite $a_{\Nis}\omega_{\CI}\colon \CI\to \RSC_{\Nis}$ 
is exact and hence commutes with homology.  
Since $\CI$ is abelian, $\mathrm{HH}_n^{(m)}(\GaSharp)$ is the
homology object of $C_\bullet^{(m)}(\GaSharp)$ in $\CI$, and therefore
\begin{align}
 h_{0,\Nis}\bigl(\mathrm{HH}_n^{(m)}(\GaSharp)\bigr)
 &\simeq
 H_n\left(a_{\Nis}\omega_{\CI}C_\bullet^{(m)}(\GaSharp)\right)\notag\\
 &\simeq H_n(\Acal_{\bullet+1,m},\partial_\bullet).
 \label{eq:HH-exactness}
\end{align}
For $n=0$, we have $\partial_1=d_0-d_1=0$
by the commutativity of $\mu$.
Hence
\[
 H_0(\Acal_{\bullet+1,m},\partial_\bullet)
 =\Acal_{1,m}\simeq\Omega^m
\]
by \autoref{thm:one-additive}.
Thus we may assume $n\geq1$.
We first show that 
\begin{equation}\label{eq:HH-theta-boundary-zero}
 \theta_{[n],m}\circ\partial_{n+1}=0.
\end{equation}
Let $F/k$ be a finitely generated field and let $E/F$ be a finite
extension.  For
\[
 \alpha
 =
 [a_0,a_1,\ldots,a_{n+1};b_1,\ldots,b_m]_{E/F}
 \in\Acal_{n+2,m}(F),
\]
put $\beta:=\bigwedge_{j=1}^m\dlog b_j$.
By \eqref{eq:HH-boundary-field} and \eqref{eq:fieldbd}, we have
\begin{align*}
\theta_{[n],m,F}\,\partial_{n+1,F}(\alpha)& =
\Tr_{E/F}\Biggl(
 a_0a_1\,
 \d a_2\wedge\cdots\wedge\d a_{n+1}\wedge\beta\\
&\qquad+
 \sum_{i=1}^{n}(-1)^i
 a_0\,
 \d a_1\wedge\cdots\wedge
 \d(a_i a_{i+1})
 \wedge\cdots\wedge
 \d a_{n+1}\wedge\beta\\
&\qquad+
 (-1)^{n+1}
 a_{n+1}a_0\,
 \d a_1\wedge\cdots\wedge\d a_n\wedge\beta
 \Biggr).
\end{align*}
Using $ \d(a_i a_{i+1})=a_i\,\d a_{i+1}+a_{i+1}\,\d a_i$, 
the terms cancel successively.  More precisely, the term
$a_0a_1\,
 \d a_2\wedge\cdots\wedge\d a_{n+1}\wedge\beta$ 
cancels with the first part of the summand for $i=1$; for
$2\leq i\leq n$, the first part of the summand for $i$ cancels with
the second part of the summand for $i-1$; and the second part of the
summand for $i=n$ cancels with $ (-1)^{n+1}a_{n+1}a_0\,\d a_1\wedge\cdots\wedge\d a_n\wedge\beta$. 
Hence
\[
 \theta_{[n],m,F}\,\partial_{n+1,F}(\alpha)=0.
\]
Since such symbols generate $\Acal_{n+2,m}(F)$ by
\autoref{thm:pres}, we obtain
\[
 \theta_{[n],m,F}\circ\partial_{n+1,F}=0
\]
for every finitely generated field $F/k$.  Hence
\eqref{eq:HH-theta-boundary-zero} holds in $\RSC_{\Nis}$.

Consequently, the restriction of $\theta_{[n],m}$ to $\Ker\left(\partial_n:\Acal_{n+1,m}\to\Acal_{n,m}\right)$ 
vanishes on the image $\img\left(\partial_{n+1}:\Acal_{n+2,m}\to\Acal_{n+1,m}\right)$.
Thus $\theta_{[n],m}$ induces a canonical morphism
\begin{equation}\label{eq:HH-induced-map}
 H_n\bigl(\Acal_{\bullet+1,m},\partial_\bullet\bigr)
 \to
 \Omega^{n+m}.
\end{equation}
Together with \eqref{eq:HH-exactness}, this gives a canonical morphism
\begin{equation}\label{eq:HH-HKR-morphism}
 h_{0,\Nis}\bigl(\mathrm{HH}_n^{(m)}(\GaSharp)\bigr)
 \to
 \Omega^{n+m}.
\end{equation}

We show that \eqref{eq:HH-HKR-morphism} is an isomorphism.
Let $F/k$ be a finitely generated field.  By
\autoref{prop:HHcomplex}, the isomorphisms
$\Theta_{\bullet+1,m,F}$ give an isomorphism of complexes
\[
 \bigl(\Acal_{\bullet+1,m}(F),\partial_\bullet\bigr)
 \simeq
 \bigoplus_{q\geq0}
 C_{\bullet,q}(B,\Omega_F^{q+m}).
\]
Therefore
\begin{align*}
 H_n\bigl(\Acal_{\bullet+1,m}(F),\partial_\bullet\bigr)
 &\simeq
 \bigoplus_{q\geq0}
 H_n\bigl(C_{\bullet,q}(B,\Omega_F^{q+m})\bigr)\\
 &\simeq
 \Omega_F^{n+m},
\end{align*}
where the second isomorphism follows from \autoref{lem:HHnorm}.
Indeed, the only nonzero summand is the one with $q=n$.

We now identify this isomorphism with the morphism induced by
$\theta_{[n],m,F}$.  In degree $n$, the decomposition
\eqref{eq:HH-Theta-B} is
\[
 \Acal_{n+1,m}(F)
 \xrightarrow{\simeq}
 \bigoplus_{S\subset[n]}
 \Omega_F^{|S|+m}\otimes\Z e_S.
\]
The summand corresponding to $S=[n]$ is
\[
 \Omega_F^{n+m}\otimes\Z e_{[n]}
 =
 \Omega_F^{n+m}\otimes\Z\varepsilon^{\otimes n}.
\]
By the definition of $\Theta_{n+1,m,F}$, the projection onto this
summand is precisely
\[
 \theta_{[n],m,F}:
 \Acal_{n+1,m}(F)\to\Omega_F^{n+m}.
\]
On the other hand, by \autoref{lem:HHnorm}, the isomorphism
\[
 H_n\bigl(C_{\bullet,n}(B,\Omega_F^{n+m})\bigr)
 \simeq
 \Omega_F^{n+m}
\]
sends the class represented by $\omega\otimes\varepsilon^{\otimes n}$ to $\omega$. 
Hence the morphism induced by
$\theta_{[n],m,F}$ is exactly the above fieldwise isomorphism
\[
 H_n\bigl(\Acal_{\bullet+1,m}(F),\partial_\bullet\bigr)
 \xrightarrow{\simeq}
 \Omega_F^{n+m}.
\]
It follows that \eqref{eq:HH-HKR-morphism} is an isomorphism on every
finitely generated field over $k$.  By \cite[Lemma~4.16]{RSY22}, it is
therefore an isomorphism in $\RSC_{\Nis}$.  Thus
\[
 h_{0,\Nis}\bigl(\mathrm{HH}_n^{(m)}(\GaSharp)\bigr)
 \simeq
 \Omega^{n+m}.
\]

Finally, since this morphism is induced by $\theta_{[n],m}$, on 
a finitely generated field $F/k$ and a finite extension $E/F$, 
the isomorphism is given by 
\[
 [a_0,\ldots,a_n;b_1,\ldots,b_m]_{E/F}
 \mapsto
 \Tr_{E/F}\left(
 a_0\,\d a_1\wedge\cdots\wedge\d a_n
 \wedge
 \dlog b_1\wedge\cdots\wedge\dlog b_m
 \right).
\]
\end{proof}

\appendix
\section{Presentations of absolute differentials}\label{app:diff}

This appendix compares two presentations of absolute K\"ahler differentials
over a field.  We first recall \cite[Lemma~4.2]{BK86} and then derive an
$F$-linear form that agrees with the presentation of
\cite[Lemma~4.1]{BE03}.  The derivation from \cite{BK86} shows in particular
that the $F$-linear form below is valid over an arbitrary field.

Let $F$ be a field and let $q\geq1$.  Put
\begin{equation}\label{eq:ATq}
 T_q(F):=F\otimes_{\Z}(F^\times)^{\otimes q}.
\end{equation}
We regard $T_q(F)$ as an $F$-vector space through its first tensor factor.
There is a natural surjection
\begin{equation}\label{eq:Adelta}
 \delta_q:T_q(F)\to\Omega_F^q,
 \qquad
 a\otimes b_1\otimes\cdots\otimes b_q
 \mapsto a\,\dlog b_1\wedge\cdots\wedge \dlog b_q.
\end{equation}

By  
\cite[Lemma~4.2]{BK86}, 
the kernel of \eqref{eq:Adelta}, as an additive group, is generated by the following two
families.  The first consists of
\begin{equation}\label{eq:ABK1}
 a\otimes b_1\otimes\cdots\otimes b_q,
 \qquad b_i=b_j\quad\text{for some }i<j.
\end{equation}
The second consists of
\begin{equation}\label{eq:ABK2}
 \sum_{\nu=1}^{m}x_\nu\otimes x_\nu\otimes\beta
 -\sum_{\mu=1}^{m'}y_\mu\otimes y_\mu\otimes\beta,
 \qquad
 \sum_{\nu=1}^{m}x_\nu=\sum_{\mu=1}^{m'}y_\mu,
\end{equation}
where $x_\nu,y_\mu\in F^\times$ and
$\beta=b_2\otimes\cdots\otimes b_q$ with $b_i\in F^\times$; when
$q=1$, the factor $\beta$ is omitted.  These are
\cite[Equations~(4.2.1) and~(4.2.2)]{BK86}.


\begin{cor}\label{cor:Aflin}
For every field $F$ and every $q\geq1$, the kernel of $\delta_q$, 
as an $F$-vector subspace of $T_q(F)$, is generated by the elements
\eqref{eq:ABK1} together with
\begin{equation}\label{eq:ACat}
 x\otimes x\otimes\beta
 +(1-x)\otimes(1-x)\otimes\beta,
 \qquad x,1-x\in F^\times,
\end{equation}
where $\beta=b_2\otimes\cdots\otimes b_q$ as above and is omitted when
$q=1$.
\end{cor}

\begin{proof}
Let $N_q(F)$ be the $F$-vector subspace generated by \eqref{eq:ABK1} and
\eqref{eq:ACat}.  Every element of \eqref{eq:ABK1} lies in
$\Ker(\delta_q)$ because two logarithmic differentials coincide.  Moreover,
\begin{align*}
 \delta_q\bigl(x\otimes x\otimes\beta
 +(1-x)\otimes(1-x)\otimes\beta\bigr)
 &=(dx+d(1-x))\wedge d\log\beta=0,
\end{align*}
where $d\log\beta$ denotes the wedge of the logarithmic differentials of
the factors of $\beta$.  Thus $N_q(F)\subset\Ker(\delta_q)$.

By \cite[Lemma~4.2]{BK86}, it remains to show that every element of
\eqref{eq:ABK2} belongs to $N_q(F)$.  Fix $\beta$ and, in
$T_q(F)/N_q(F)$, write
\begin{equation}\label{eq:Ae}
 e_\beta(a):=a\otimes a\otimes\beta
 \qquad(a\in F^\times).
\end{equation}
If $a,b,a+b$ are nonzero, apply \eqref{eq:ACat} with
$x=a/(a+b)$ and multiply the first tensor factor by $a+b$.  We obtain
\[
 a\otimes\frac{a}{a+b}\otimes\beta
 +b\otimes\frac{b}{a+b}\otimes\beta=0.
\]
Since $F^\times$ is written additively inside the tensor product,
$a/(a+b)=a-(a+b)$ and $b/(a+b)=b-(a+b)$.  Hence
\begin{equation}\label{eq:Aadd}
 e_\beta(a)+e_\beta(b)=e_\beta(a+b)
 \qquad(a,b,a+b\neq0).
\end{equation}

If $a+b=0$ and $F\neq\mathbb F_2$, choose
$c\in F^\times$ with $c\neq-a$.  Applying \eqref{eq:Aadd} to $(a,c)$
and to $(-a,a+c)$ gives
\[
 e_\beta(a)+e_\beta(c)=e_\beta(a+c),
 \qquad
 e_\beta(-a)+e_\beta(a+c)=e_\beta(c),
\]
so $e_\beta(-a)=-e_\beta(a)$.  With $e_\beta(0):=0$, the map
$e_\beta:F\to T_q(F)/N_q(F)$ is therefore additive.  When
$F=\mathbb F_2$, one has $F^\times=\{1\}$ and $T_q(F)=0$ for $q\geq1$.

Now the equality of sums in \eqref{eq:ABK2} gives
\[
 \sum_{\nu=1}^{m}e_\beta(x_\nu)
 =e_\beta\left(\sum_{\nu=1}^{m}x_\nu\right)
 =e_\beta\left(\sum_{\mu=1}^{m'}y_\mu\right)
 =\sum_{\mu=1}^{m'}e_\beta(y_\mu).
\]
Thus every element of \eqref{eq:ABK2} lies in $N_q(F)$.  Together with
\cite[Lemma~4.2]{BK86}, this proves $N_q(F)=\Ker(\delta_q)$.
\end{proof}

The relation \eqref{eq:ACat} is also directly a special case of
\eqref{eq:ABK2}: take $x_1=x$, $x_2=1-x$, and $y_1=1$; the term
$1\otimes1\otimes\beta$ is zero because the second tensor factor is the
identity of $F^\times$.  Thus the passage from \eqref{eq:ABK2} to
\eqref{eq:ACat} amounts precisely to using $F$-linear combinations to
reduce finite additivity to the normalized identity $x+(1-x)=1$.

Bloch and Esnault give the corresponding $F$-linear presentation in
\cite[Lemma~4.1]{BE03}, using
$F\otimes_{\Z}\bigwedge^qF^\times$.  Passing from the exterior power to
the tensor power introduces the repeated-entry relations \eqref{eq:ABK1},
and the remaining relation is exactly \eqref{eq:ACat}.  Hence
\autoref{cor:Aflin} identifies this presentation with that of
\cite[Lemma~4.2]{BK86} in the field case.  The tensor form
\eqref{eq:ACat} is also used in \cite[p.~497, (11)]{Hir14}.


\def\cprime{$'$}
\providecommand{\bysame}{\leavevmode\hbox to3em{\hrulefill}\thinspace}

\end{document}